\documentclass[preprint,12pt]{elsarticle}

\usepackage{amssymb}
\usepackage{amsmath}

\usepackage{subcaption}
\usepackage{algorithm}
\usepackage{algorithmic}
\DeclareUnicodeCharacter{200B}{}
\usepackage{booktabs}
\usepackage{multirow}
\usepackage{float}
\usepackage[export]{adjustbox}
\usepackage[left=2cm,right=2cm,top=2.5cm,bottom=2.5cm]{geometry}
 \usepackage{lineno}

\journal{Inverse Problems and Imaging}

\begin{document}

\begin{frontmatter}



\title{DLMMPR:	Deep Learning-based Measurement Matrix for Phase Retrieval}


\author[label1]{Jing Liu} 
\author[label1]{Bing Guo}
\author[label1]{Ren Zhu}
\affiliation[label1]{organization={School of Mathematics and Statistics},
            addressline={Jishou University}, 
            city={Jishou},
            postcode={416000}, 
            state={Hunan Province},
            country={China}}

\begin{abstract}
  This paper pioneers the integration of learning optimization into measurement matrix design for phase retrieval. We introduce the Deep Learning-based Measurement Matrix for Phase Retrieval (DLMMPR) algorithm, which parameterizes the measurement matrix within an end-to-end deep learning architecture. Synergistically augmented with subgradient descent and proximal mapping modules for robust recovery, DLMMPR’s efficacy is decisively confirmed through comprehensive empirical validation across diverse noise regimes. Benchmarked against DeepMMSE and PrComplex, our method yields substantial gains in PSNR and SSIM, underscoring its superiority.
 
\end{abstract}


\end{frontmatter}



\section{Introduction}

The characterization of a monochromatic coherent light field in optics relies on its two-dimensional complex amplitude function, which encodes both amplitude and phase. However, existing optical imaging devices can only acquire the intensity  information, but face significant hurdles in directly capturing phase information. This fundamental challenge necessitates the development of phase retrieval technology, a field devoted to recovering this lost or degraded phase information. Through the measurement of light field intensity distributions, phase retrieval employs numerical algorithms to reconstruct the complete wave phase. Positioned at the nexus of computational mathematics and data science, this powerful technique finds widespread application in diverse areas such as diffraction imaging \cite{gerchberg1994practical}, biomedical imaging \cite{kong2023phase}, astronomical observations \cite{fienup1987phase}, electron microscopy \cite{zhang2021physics}, and coherent inverse scattering \cite{metzler2017coherent}.

 The core of the phase retrieval problem lies in reversely deducing the complete complex amplitude characteristics  from light intensity information. Specifically, it refers to reconstructing  an unknown signal from a given set of measured values of amplitude or intensity.
Let $x \in \mathbb{C}^N$represent the signal to be retrieved. The measured value of the signal, denoted as  $\langle a_i,x\rangle$, is obtained via a known sampling vector $a_i \in \mathbb{C}^N$.

Define the sampling matrix $\boldsymbol{A}^T = [a_1, a_2, \ldots, a_m] \in \mathbb{C}^{N \times M}$ and the magnitude vector $\boldsymbol{y}^T = [y_1, y_2, \ldots, y_m] \in \mathbb{C}^{1 \times M}$, where 
$w$ denotes noise. Then, the phase retrieval problem can be  expressed as
\begin{equation}
	\boldsymbol{y} = \vert \boldsymbol{A x} \vert + w.
	\label{PR} 
\end{equation}

To address the issue of phase loss, researchers in the 1970s proposed a framework based on alternating projection technology to achieve phase retrieval.  Gerhberg and Saxton \cite{gerhberg1972practical} put forward the earliest GS (Gerchberg-Saxton) algorithm, which has since become a classic algorithm for solving phase retrieval problems.  Fienup \cite{fienup1984reconstruction} extended the GS algorithm and proposed the HIO algorithm (Hybrid Input-Output algorithm).
These methods are sensitive to noise and cannot converge quickly. To tackle these problems, scholars have developed several improved algorithms. Examples include the Hybrid Projection Reflection (HPR) algorithm \cite{bauschke2003hybrid}, the Relaxed Averaged Alternating Reflections (RAAR) algorithm \cite{luke2004relaxed}, and the Difference Map (DM) algorithm \cite{elser2003phase}. These improved algorithms are designed to handle different types of specific problems and further enhance the convergence of the algorithms.

 Many different theories have been attempted to be integrated into phase retrieval. Compressive sensing phase retrieval was proposed in \cite{moravec2007compressive}, which introduced compressibility constraints based on sparse bases. It provided a  technically feasible solution to break through the limitations of traditional phase retrieval methods.  Newton et al. \cite{newton2012compressed}  proposed a compressive sensing HIO algorithm (CSHIO) that adopts reweighted 
$\ell_1$-norm minimization. It can be used for phase retrieval  in non-convex scenarios. Both of these algorithms are designed for real-valued distributions.  Based on the Alternating Direction Method of Multipliers (ADMM), Chang et al. \cite{chang2018total} presented  an efficient numerical algorithm with total variation regularization. However, the total variation weight parameter 
$\lambda$ and the penalty parameter need to be adjusted manually, and its sensitivity and computational complexity limit its practical application.

Deep learning methodologies, especially when integrated into the Plug-and-Play Alternating Direction Method of Multipliers (PnP-ADMM) framework \cite{venkatakrishnan2013plug,chan2016plug}, have revolutionized phase retrieval applications. Kappeler \cite{kappeler2017ptychnet} pioneered PtychNet, a Convolutional Neural Network (CNN) that dramatically curtails the computational burden of iterative phase retrieval algorithms. Christopher et al. \cite{metzler2018prdeep}  further refined this by proposing prDeep, which masterfully fuses the Denoising Regularization (RED) framework with phase retrieval, embedding a pre-trained DnCNN denoiser within the PnP structure. This synergistic approach not only drastically slashed iteration times but also achieved exceptional recovery performance. However, a persistent limitation of existing PnP methods is their reliance on manual parameter adjustment, rendering them susceptible to yielding suboptimal results when faced with significant variations in imaging conditions or scene content. To surmount this hurdle, Wei et al. \cite{wei2020tuning} introduced a parameter-free PnP proximal algorithm capable of autonomously determining key parameters like penalty coefficients, denoising strength, and stopping thresholds. Its core innovation lies in a dedicated strategy network that orchestrates automatic parameter discovery through a sophisticated blend of model-free and model-based deep reinforcement learning. Most recently, PRISTA-Net \cite{liu2023prista} has emerged as a novel Deep Unfolding Network (DUN), built upon the foundation of the first-order Iterative Shrinkage-Thresholding Algorithm (ISTA). This network leverages learnable nonlinear transformations to elegantly solve the proximal mapping subproblem. Crucially, all intrinsic parameters—including nonlinear transformations, thresholding values, and step sizes—are optimized end-to-end, completely eliminating the need for manual intervention.

Most of the  phase retrieval research focuses on the optimization of reconstruction algorithms, such as iterative strategies, denoising regularization, and deep learning network design. However, as the front-end of information acquisition, the measurement matrix is usually constructed based on physical models and remains dominated by fixed designs. For example, Coded Diffraction Patterns (CDPs) generate matrix 
$A$ through the combination of random masks and Fourier transforms, and its structure is not optimized for specific image distributions or imaging scenarios, leading to the following inherent drawbacks.

\textbf{Low Information Acquisition Efficiency}~
Fixed matrices struggle to match the sparse structures and phase characteristics of different data types, such as natural images and medical images. Under low sampling rates, they tend to lose key information, which directly limits the room for improving reconstruction accuracy.

\textbf{Insufficient Noise Robustness}~
The design of fixed matrices does not take noise distribution characteristics into account. In high-noise environments, the signal-to-noise ratio (SNR) of measured data is further reduced, which increases the difficulty of subsequent reconstruction.

In response to the critical challenge of measurement matrix optimization, Abolghasemi et al. \cite{measurement} illuminated the pivotal role of mutual coherence between the measurement matrix and its representation counterpart in influencing signal reconstruction. They advocated for the gradient descent methodology,which aimed at diminishing this coherence by minimizing the absolute values of the off-diagonal elements, thereby offering foundational theoretical guidance for matrix design. Chen et al. \cite{2016Measurement} transcended the constraints of static matrix design by introducing a novel optimization strategy predicated on genetic algorithms. This  approach not only accounts for the practicalities of physical observation and target attributes, but also formulates the measurement matrix as a randomly chosen partial unitary matrix.
 Concurrently, 
  relevant research \cite{2023Measurement} adeptly leveraged prior target imaging information.
   and yielded both optimal measurement matrices and high-resolution imaging outcomes, while achieving a substantial enhancement in the downsampling rate of radar data.

The foregoing review confirms that optimizing the measurement matrix is central to advancing sparse imaging system performance. Consequently, this paper reorients focus towards measurement matrix optimization in phase retrieval, circumventing the limitations of fixed designs. Employing a learning-based strategy for adaptive engagement with intricate image distributions and noise conditions, we aim to harness a measurement matrix performance leap. This endeavor seeks to provide novel solutions for high-fidelity, efficient, and robust phase retrieval, driving collaborative optimization across the entire acquisition-to-reconstruction continuum.

The main contributions of this work are the proposal of the Deep Learning-based Measurement Matrix for Phase Retrieval (DLMMPR) algorithm and its associated reconstruction methodology. DLMMPR articulates a sophisticated, learnable deep neural network, architected upon the structural blueprint of ISTA update steps. Critically, the measurement matrix is dynamically learned from prior information, eschewing manual design or fixed configurations, thereby substantially advancing phase retrieval performance. This paradigm shift endows DLMMPR with the dual advantages of swift, precise reconstruction and commendable interpretability. Extensive empirical validation across natural and non-natural images decisively confirms DLMMPR’s superiority over current state-of-the-art methods, while preserving a favorable computational complexity.

\section{Background}

Learning optimization denotes the enhancement or reconstruction of classical optimization algorithms by leveraging deep learning technologies. Its primary goal is to yield more efficient optimization methods that excel in accuracy, convergence speed, and computational cost compared to traditional techniques. This powerful paradigm has demonstrated considerable promise in numerous optimization problems and fields, such as compressive sensing and inverse problems \cite{chen2018theoretical}, neural network training \cite{andrychowicz2016learning}, combinatorial optimization \cite{khalil2017learning}, computer vision \cite{zhang2019real}, and software engineering \cite{agrawal2020better}. The core advantages it confers are follows.

\textbf{Accuracy Improvement}~
It uses data-driven approaches to learn optimal iteration rules, avoiding errors in traditional algorithms caused by prior assumptions.

\textbf{Computational Lightweight}~
It leverages the parameterized mapping of neural networks to replace the complex iterative steps of traditional algorithms.

\textbf{Rapid Convergence}~
It adaptively learns the structural characteristics of problems, reducing the number of iterations. In some scenarios, the convergence speed can be increased by more than 10 times.

The genesis of learning optimization can be traced to 2010, with the pioneering contributions of Gregor et al. \cite{gregor2010learning} and their introduction of the LCoD and LISTA algorithms. This marked a watershed moment, representing the first instance of supervised learning being integrated into the critical phase of generating initial solutions for optimization, definitively proving that learned approaches could dramatically accelerate the convergence of established algorithms. LISTA, in particular, ingeniously unrolled the iterative ISTA process into a neural network, enabling supervised learning to optimize thresholds and step sizes, which resulted in an impressive threefold speedup in compressive sensing reconstruction. Simultaneously, their LCoD algorithm demonstrated, for the first time, that learned initial solutions could significantly truncate the iteration count of traditional CoD algorithms, thereby forging a new path for learning-based initialization. 

Andrychowicz et al. \cite{andrychowicz2016learning} parameterized optimization rules using LSTM, enabling the optimizer to independently learn problem characteristics and achieve adaptive gradient updates. Its core formula is $\theta_{t+1} = \theta_t + g_t(\nabla f(\theta_t), \phi)$, where $g_t$ is generated by LSTM and 
$\phi$ represents the optimizer parameters. This realizes adaptive optimization in target task classes and improves the generalization ability of resolution and style.
Chen et al. \cite{chen2018theoretical}developed the LISTA model based on ISTA. By transforming the iterative process into a trainable neural network structure, it achieved more efficient convergence performance than traditional iterative algorithms in sparse vector recovery tasks. They also proposed a key partial weight coupling structure, i.e., $\boldsymbol{W}^k_2 = \boldsymbol{I} - \boldsymbol{W}^k_1\boldsymbol{A}$. Under this structure, LISTA can achieve linear convergence — a convergence speed significantly superior to the sublinear convergence of ISTA/FISTA under general conditions.

 DualNet \cite{DualPRNet2022} integrated a new type of Deep Shrinkage Network (DSN) into a supervised dual-frame learning framework, and proposed a deep shrinkage dual-frame network for constructing a deep unfolded phase retrieval network architecture. This enables the network to learn instance-adaptive and spatially varying threshold functions, addressing the drawback of traditional threshold functions—where manually designed thresholds are used to filter frame coefficients, lacking adaptability and thus limiting the final reconstruction quality.
All parameters in the PRISTA-Net \cite{liu2023prista} framework (including nonlinear transformations, threshold parameters, and step sizes) are obtained through end-to-end learning, abandoning the parameter setting mode that relies on manual experience in traditional methods. This design allows the network to independently learn the optimal parameter configuration from data.
 Yang et al. \cite{YANG2023HIONet} proposed a deep unfolded network HIONet based on deep priors and its enhanced version HIONet$^+$ to solve the phase retrieval problem. Derived from the unfolding of the classic HIO algorithm, HIONet$^+$ replaces the projection operator with a trainable PPModule to adaptively learn the inherent priors of data, integrating the representational capability of deep networks with the interpretability and universality of traditional algorithms. Additionally, HIONet$^+$ introduces cross-block feature fusion to improve reconstruction performance.

To effectively tackle the regularization challenges in phase retrieval, we begin by examining the compressive sensing problem and its solution paradigms, referencing the ISTA algorithm’s \cite{gregor2010learning} approach to regularization for compressive sensing. The mathematical model governing compressive sensing can be precisely defined as
\begin{equation}
	\boldsymbol{y} = \boldsymbol{A x}  + w,
\end{equation}
where $\boldsymbol{y} \in \mathbb{R}^m$ denotes the observation vector, $\boldsymbol{A} \in \mathbb{R}^{m \times n}$  is the measurement matrix, $\boldsymbol{x} \in \mathbb{R}^n$ represents the sparse signal to be recovered, $w$ stands for noise.
To solve the above problem, the LASSO optimization framework is usually adopted
\begin{equation}
	\underset{x}{\text{minimize}} \frac{1}{2}\| \boldsymbol{A x} - \boldsymbol{y}\| _{2}^{2} + \lambda\| \boldsymbol{x}\| _{1}.
\end{equation}

\textbf{The traditional ISTA algorithm}~~According to the literature \cite{beck2009fast}, the classic ISTA algorithm solves this problem through iterative updates. The gradient iteration formula can be regarded as the proximal regularization of the linearized function of $f$ at $\boldsymbol{x}^{k-1}$, and its equivalent form is
\begin{equation}
	\boldsymbol{x}^k = \arg\min_{\boldsymbol{x}} \left\{ f(\boldsymbol{x}^{k-1}) + \langle \boldsymbol{x} - \boldsymbol{x}^{k-1}, \nabla f(\boldsymbol{x}^{k-1}) \rangle + \frac{1}{2t^k} \|\boldsymbol{x} - \boldsymbol{x}^{k-1}\|^2 + \lambda \|\boldsymbol{x}\|_1 \right\}.
\end{equation}
After ignoring the constant term, it can be rewritten as
\begin{equation}
	\boldsymbol{x}^k = \arg\min_{\boldsymbol{x}} \left\{ \frac{1}{2t^k} \|\boldsymbol{x} - (\boldsymbol{x}^{k-1} - t^k \nabla f(\boldsymbol{x}^{k-1}))\|^2 + \lambda \|\boldsymbol{x}\|_1 \right\}.
\end{equation}
Finally, thresholding is performed
\begin{equation}
	\boldsymbol{x}^k= \eta_{\phi^k}\left(\boldsymbol{x}^{k} - t^k \boldsymbol{A}^T (\boldsymbol{Ax}^k - \boldsymbol{y})\right), \quad k = 0, 1, 2, \ldots
\end{equation}

\begin{algorithm}[H]
	\caption{ISTA(Iterative Shrinkage and Thresholding Algorithm)\cite{gregor2010learning}}
	\label{ISTA algorithm}
	\begin{algorithmic}[1]
		\REQUIRE  Let $\boldsymbol{x}_0 $
		be the initial point, set $k :=0$, and let 
		$L$ be the maximum eigenvalue of $\boldsymbol{A}^T\boldsymbol{A}$
		\FOR{$ k=1,2,3,\ldots ,$maxiter}
		\STATE $\boldsymbol{r}^{k}=\boldsymbol{x}^{k-1}-t^k\boldsymbol{A}^{T}\left(\boldsymbol{A x}^{k-1} - \boldsymbol{y}\right)$.
		\STATE $\boldsymbol{x}^{k}=\eta_{\phi^k}\left(\boldsymbol{r}^{k}\right)$.
		\ENDFOR
		\STATE Output the current point $\boldsymbol{x}^{maxiter}$ as the approximate optimal solution..
	\end{algorithmic}
\end{algorithm}

\textbf{ISTA-Net Algorithm}~~According to ISTA-Net \cite{zhang2018ista}, prior sparsity is imposed on the $\ell_1$
norm of $\boldsymbol{x}$  in a certain transform domain, then regularization term is $ \|\boldsymbol{\varPsi}\boldsymbol{x}\|_1$. Let the transform domain $\boldsymbol{\varPsi} = \mathcal{F}$  where  $\mathcal{F}$ is a wavelet transform matrix with orthogonality, Thus, the optimization model becomes

\begin{equation}
	\underset{x}{\text{minimize}} \frac{1}{2}\| \boldsymbol{A x} - \boldsymbol{y}\| _{2}^{2} + \lambda\| \boldsymbol{\varPsi}\boldsymbol{x}\| _{1}.
\end{equation}
Therefore, the solution to the above problem is equivalent to a two-step iterative proximal gradient descent algorithm:
\begin{align}
	\boldsymbol{r}^k &= \boldsymbol{x}^{k-1} - t^k \nabla f(\boldsymbol{x}^{k-1}) \label{eq:gradient_step},\\
	\begin{split}
		\boldsymbol{x}^k &= \boldsymbol{\text{Prox}}_{\phi^k \|\mathcal{F}(\cdot)\|_1} (\boldsymbol{r}^k)\\
		&= \arg\min_{\boldsymbol{x}} \frac{1}{2} \|\mathcal{F}(\boldsymbol{x}) - \mathcal{F}(\boldsymbol{r}^k)\|_2^2 + \eta^{k} \|\mathcal{F}(\boldsymbol{x})\|_1 \\
		&=\tilde{\mathcal{F}} \text{soft}(\mathcal{F}( \boldsymbol{r}^k), \eta^{k}).
	\end{split} \label{eq:proximal_step}
\end{align}

Where the wavelet transform matrix  $\mathcal{F}$ is invertible, and its left inverse is denoted as $\tilde{\mathcal{F}}$. $\operatorname{soft}(z, \tau) = \operatorname{sign}(z) \max\{|z| - \tau, 0\}$ is the soft-thresholding function, 
$\eta^{k}$ represents the soft threshold, and  $t^k$ is the step size in the traditional composite algorithm. All learnable parameters in ISTA-Net are denoted as
$\boldsymbol{\Theta} = \left\{ t^{k}, \eta^{k}, \mathcal{F}^{k}, \tilde{\mathcal{F}}^{k} \right\}_{k=1}^{K}$, For more detailed explanations, refer to \cite{zhang2018ista}.

\subsection{The proposed algorithm}
The phase retrieval process mainly consists of two key parts: signal measurement and phase recovery. The focus of this paper is on the first part, i.e., attempting to improve the measurement process to enhance phase retrieval performance. The main flow of the proposed DLMMPR algorithm is shown in Figure \ref{fig:network}. Its core is to transform the fixed steps of traditional iterative algorithms into a trainable deep network structure, replacing manually designed measurement matrices with data-driven optimization.

Generally, the generalized phase retrieval  problem based on regularization priors can be expressed as the following optimization problem
\begin{equation}
	\arg\min_{x\in\mathcal{X}} f(\boldsymbol{x}) + \lambda g(\boldsymbol{x}).
\end{equation}

Here, $f(\boldsymbol{x}) = \frac{1}{2} \lVert |\boldsymbol{A} \boldsymbol{x}| - \boldsymbol{y} \rVert_2^2$ is the data fidelity term, $g(\boldsymbol{x})$ serving as the regularization term, represents the prior information of the data, and the regularization parameter $\lambda$ denotes the penalty degree imposed on $g(\boldsymbol{x})$ First, based on the Iterative Shrinkage-Thresholding Algorithm (ISTA\cite{zhang2018ista}), we improve the phase retrieval algorithm. ISTA is used to impose prior sparsity on the $\ell_1$-norm of $\boldsymbol{x}$ in a specific transform domain, i.e., $\ g(\boldsymbol{x}) = \|\boldsymbol{\varPsi}\boldsymbol{x}\|_1$,where
$\boldsymbol{\varPsi}\boldsymbol{x}$ represents the transform coefficients of $\boldsymbol{x}$ with respect to a certain transform  $\boldsymbol{\varPsi}$. By incorporating the sparsity constraint as a regularization term into the smooth non-convex objective function, a regularized objective function is obtained. Therefore, the regularized phase retrieval model that relies on amplitude measurements and uses the  $\ell_1$-norm to characterize sparse signals is as follows
\begin{equation}
	\label{model}  \frac{1}{2}\||\boldsymbol{A}\boldsymbol{x}|-\boldsymbol{y}\|_2^2+\lambda \|\boldsymbol{\varPsi}\boldsymbol{x}\|_1.
\end{equation}

\subsection{Signal measurement}

\begin{figure}[!t]
	\centering
	\renewcommand{\figurename}{Figure}
	\includegraphics[width=1\textwidth]{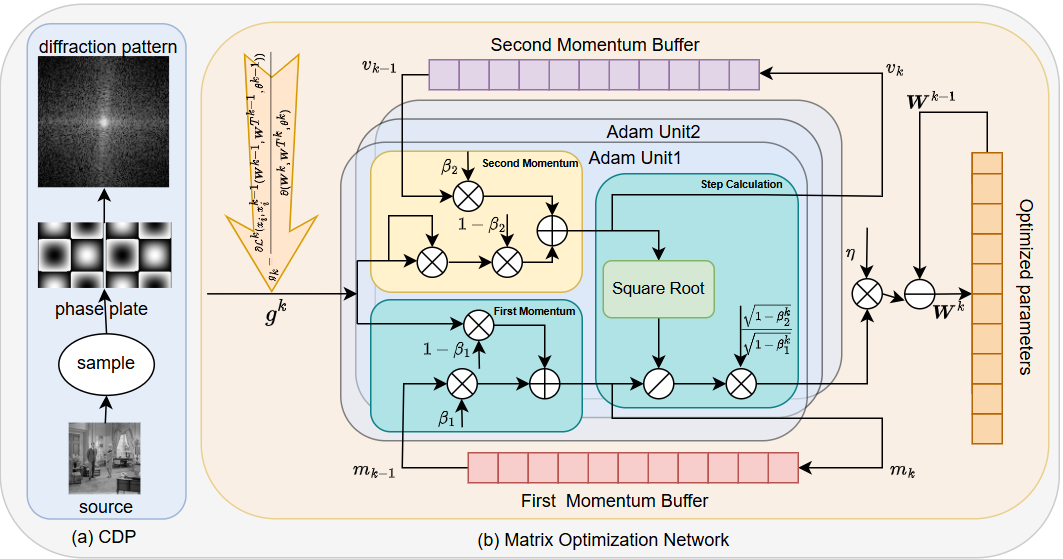}
	\caption{Coded Diffraction Patterns and Matrix Optimization Network.}
	\label{fig:network-adam} 
\end{figure}

 The fidelity of phase reconstruction in phase retrieval problems is fundamentally dictated by the design of the measurement method. While Fourier measurement and Coded Diffraction Phase (CDP) measurement are the prevailing techniques for most phase retrieval applications, this paper focuses on enhancing the CDP framework (illustrated in Figure  \ref{fig:network-adam} (a)). We address the crucial problem of learning an optimal measurement matrix within the CDP paradigm, thereby proposing significant improvements over traditional phase retrieval models’ measurement strategies.
  Specifically, $\boldsymbol{A}$ and its conjugate transpose  $\boldsymbol{A}^\mathrm{H}$ are added as learning parameters, and the specific update process of the learning parameters is shown in Fig. \ref{fig:network-adam} (b).

\textbf{Traditional Measurement Algorithm}~ The Coded Diffraction Pattern (CDP) technique \cite{candes2015phase} harnesses the power of Spatial Light Modulators (SLMs) to strategically expand the frequency spectrum of the target, significantly simplifying subsequent reconstruction efforts. Within this framework, the signal is modulated before it diffracts. Common strategies for designing these crucial scheduling waveforms include

\textbf{\textbullet~ Phase Mask Implementation:}~Introducing a phase mask in proximity to the sample to impart targeted phase shifts.

\textbf{\textbullet~ Optical Grating Modulation:}~Employing an optical grating to sculpt the illuminating beam, thereby precisely controlling the phase or amplitude characteristics of the light field.

\textbf{\textbullet~ Stochastic Modulation:}~Utilizing randomly generated waveforms, such as complex Gaussian or binary distributions, to introduce stochastic variations.

\textbf{\textbullet~ Ptychographic Scanning:}~Achieving modulation by systematically scanning an illuminating spot across an extended sample, enabling the capture of overlapping diffracted fields.

Considering a scheduled waveform denoted as  $d[t]$ and an original signal ${x[t]}$
, the observed diffraction pattern, known as the Coded Diffraction Pattern (CDP), represents the spectral content of ${x[t]}$ as modulated by the code  $d[t]$. This interaction is mathematically described as follows
\begin{equation}
	y_k = \left| \sum_{t=0}^{n-1} x[t] \bar{d}[t] e^{-i 2\pi \omega_k t} \right|^2, \quad \omega_k \in \Omega.
\end{equation}

Let ${d_\ell[t]}$ denote the  ${\ell}$-th scheduling waveform, where $k$
corresponds to different frequency components of the Discrete Fourier Transform (DFT), and $L$ is the number of scheduling times. After  $L$ modulations, the DFT modulus of the signal $\boldsymbol{x}$ forms the "encoded" diffraction pattern as follows
\begin{equation}
	y_{\ell,k} = \left| \sum_{t=0}^{n-1} x[t] \bar{d}_\ell[t] e^{-i 2\pi k t / n} \right|^2, 
	\quad
	\begin{aligned}
		&0 \leq k \leq n - 1, \\ 
		&1 \leq \ell \leq L. 
	\end{aligned}
\end{equation}

Expressed in matrix notation, let   $\boldsymbol{D}_\ell$ be a diagonal matrix whose diagonal elements are the modulation patterns $d_\ell [t]$, and $\boldsymbol{f}_k^*$ be the row vector of the DFT matrix. Then the measured value can be expressed as
\begin{equation}
	y_{\ell,k} = \left| \boldsymbol{f}_k^* \boldsymbol{D}_\ell \boldsymbol{x} \right|^2.
\end{equation}

The CDP encoding process intrinsically embeds additional structural constraints within the measured data. This distinguishes it from Fourier measurement, which provides only single Fourier modulus information, making CDP inherently more adaptable to complex environments. Its superior adaptability is particularly evident in scenarios involving weak signals, noisy conditions, or the need for higher measurement efficiency. Consequently, CDP drives phase retrieval’s progression from empirical algorithms to a domain grounded in rigorous mathematical theory, fostering a more complete theoretical system for the field.

The measurement model for phase retrieval is described by Model (\ref{PR}), where the measurement matrix $\boldsymbol{A}$ is constructed in the form of multi-channel Fourier modulation as
\begin{equation*}
 \boldsymbol{A} = \left[ (\boldsymbol{F}\boldsymbol{D}_1)^T, (\boldsymbol{F}\boldsymbol{D}_2)^T, \cdots, (\boldsymbol{F}\boldsymbol{D}_J)^T \right]^T,
 \end{equation*}
 here,    $\boldsymbol{F}$ is a two-dimensional Fourier transform matrix, and $\{\boldsymbol{D}_j\}_{j=1}^J$ are diagonal illumination mask matrices simulating Spatial Light Modulators (SLMs). Their diagonal elements are sampled from the unit circle in the complex plane (with a modulus of 1 and uniformly distributed arguments), ensuring the randomness and coverage of the modulation phases. In actual measurements, an additive noise term $\boldsymbol{w}$ is introduced due to the influence of hardware noise. Based on the verification of the robustness of multi-measurement phase retrieval in reference \cite{metzler2018prdeep}, this paper uses $J=4$ to simulate four independent measurements of a pure-phase SLM, where 
$\boldsymbol{D}_1, \boldsymbol{D}_2, \boldsymbol{D}_3, \boldsymbol{D}_4$ correspond to the diagonal matrices constructed as described above.

\textbf{Measurement Matrix Optimization}~~ It is worth noting that in traditional CDPs, $\boldsymbol{A}$ is often used as a fixed matrix. Inspired by the idea of LISTA, this paper explores the possibility of embedding $\boldsymbol{A}$
and its conjugate transpose $\boldsymbol{A}^\mathrm{H}$ as learnable parameters into an end-to-end framework to meet the phase retrieval requirements in complex scenarios.

Before the start of training, we use the model (see \ref{PR}) to calculate the corresponding phase retrieval measurements $y_i$. The measurement matrix $\boldsymbol{A}$ here remains the traditional CDP measurement, which is obtained directly using the fast Fourier function. In this way, we construct our training dataset  $\mathcal{D}_{train} = \left\{ (y_i, x_i) \right\}_{i=1}^{\mathcal{N}_s}$ 
using a large number of image patches and their corresponding phase retrieval measurements. Then, we define a learnable measurement matrix, perform learning and optimization on the 2D Fourier transform matrix $\boldsymbol{F}$, and denote the optimized Fourier transform matrix as $\boldsymbol{T}$, which is expressed as
\begin{equation}
	\boldsymbol{W} = \left[ (\boldsymbol{T}\boldsymbol{D}_1)^T, (\boldsymbol{T}\boldsymbol{D}_2)^T,(\boldsymbol{T}\boldsymbol{D}_3)^T , (\boldsymbol{T}\boldsymbol{D}_4)^T \right]^T,
\end{equation}
denoting the learned matrix $\boldsymbol{A}$ as $\boldsymbol{W}$
and the learned matrix  $\boldsymbol{A}^T$ as  $\boldsymbol{W}^T$
respectively. The improved measurement model is expressed as
\begin{equation}
	\boldsymbol{y} = \vert \boldsymbol{W} \boldsymbol{x} \vert + \boldsymbol{w}.
	\label{measurement} 
\end{equation}

All phase retrieval iterations are performed contingent upon the model described in Figure \ref{measurement}. Our methodology integrates learned measurement matrices, 
 $\boldsymbol{W}$,and $\boldsymbol{W}^T$, along with a set of other learnable parameters represented by $\theta$, encompassing regularization weights, iteration step sizes, and auxiliary configurations. Parameter updates are achieved through an end-to-end learning framework for adaptive optimization. Using the phase retrieval measurement  $y_i$as network input and the original image patch $x_i$
as the target output, the multi-stage iterative process generates a reconstruction 
 $x_i^{\mathcal{N}_s}$. Parameter adjustment is driven by minimizing the difference between the reconstruction  $x_i^{\mathcal{N}_s}$ and the ground truth  $x_i$. This reconstruction error is quantified by a well-posed loss function, defined as follows
\begin{equation}
	\mathcal{L}^k = \frac{1}{\mathcal{N}_s \mathcal{N}} \sum_{i=1}^{\mathcal{N}_s} \|x_i^k (\boldsymbol{W}^k, \boldsymbol{W}^{T^k}, \theta^k) - x_i\|_2^2
	\label{eq:loss}.
\end{equation}

The gradients of the loss function with respect to all learnable parameters including  $\boldsymbol{W}$ and $\boldsymbol{W}^T$
are calculated using the backpropagation algorithm. An Adam optimizer is designed, and the parameter update process of Adam is detailed in Algorithm \ref{Adam algorithm}. In each iteration, the measurement matrices 
$\boldsymbol{W}$ and $\boldsymbol{W}^T$ are first updated to minimize the measurement error; then, the auxiliary parameters $\theta$
are optimized based on the current matrix parameters to enhance the effectiveness of prior constraints. By simultaneously calculating the first moment 
$\beta_1$ and the second moment $\beta_2$ of the gradients, a unique learning rate is computed for each parameter. The initial hyperparameters are set as 
$\beta_1=0.9$,$\beta_2=0.999$, and step size $\eta=0.01$. This ultimately achieves a global optimal balance of all parameters during the convergence process, enabling adaptive adjustment of parameter values.

In the realm of measurement matrix learning, our end-to-end framework ushers in an era of data-driven, adaptive learning. It liberates us from the constraints of manually defined feature extraction rules, enabling the automatic discovery of optimal feature interdependencies and measurement logic directly from the data. This sophisticated process effectively constructs a dynamically adaptive pathway for data processing. Crucially, by streamlining away multi-stage intermediate conversions, the framework inherently curtails information loss. When synergized with the accelerated parameter updates of the Adam optimizer, this leads to a substantial leap in reconstruction accuracy. For enhanced noise robustness, the end-to-end framework is engineered to intrinsically incorporate noise suppression capabilities alongside the learning of the reconstruction task via joint training. It forges a direct, end-to-end mapping from noisy measurements to clean images, thereby precluding the detrimental accumulation and amplification of noise typically seen in sequential processing. Furthermore, the Adam optimizer guarantees stable model convergence even under challenging noise conditions, markedly improving the phase retrieval process’s resilience. Subsequent experimental validation robustly confirms that our method adeptly overcomes the technical limitations of conventional fixed matrix approaches, delivering a more efficient and robust solution for phase retrieval endeavors.
\begin{algorithm}[H]
	\caption{Adam algorithm}
	\label{Adam algorithm}
	\begin{algorithmic}[1]
		\REQUIRE  Exponential decay rate of moment estimation $\beta_1=0.9$,$\beta_2=0.999$
		,Step Size $\eta=0.01$,$\epsilon=10^{-8}$
		, Initialize the first-moment vector and the second-moment vector $m_0 \leftarrow 0$ and $v_0 \leftarrow 0$
		\REQUIRE  $x_i$,$ x_i^{k-1} (\boldsymbol{W}^{k-1}, \boldsymbol{W}^{T^{k-1}})$
		\FOR{$ k=1,2,3,\ldots ,$K}
		\STATE $g_k = \frac{\partial \mathcal{L}^k(x_i,x_i^{k-1} (\boldsymbol{W}^{k-1}, \boldsymbol{W}^{T^{k-1}}, \theta^{k-1}))}{\partial (\boldsymbol{W}^k, \boldsymbol{W}^{T^k}, \theta^k)}$,
		\STATE $m_k = \beta_1 \cdot m_{k-1} + (1 - \beta_1) \cdot g_k$.
		\STATE $v_k = \beta_2 \cdot v_{k-1} + (1 - \beta_2) \cdot g_k^2$,
		\STATE$\hat{m}_k = \frac{m_k}{1 - \beta_1^k}, \quad \hat{v}_k = \frac{v_k}{1 - \beta_2^k}$,
		\STATE $\boldsymbol{W}^k=\boldsymbol{W}^{k-1} - \eta \cdot \frac{\hat{m}_k}{\sqrt{\hat{v}_k} + \epsilon},\boldsymbol{W}^{T^k}=\boldsymbol{W}^{T^k-1} - \eta \cdot \frac{\hat{m}_t}{\sqrt{\hat{v}_t} + \epsilon}$.
		\ENDFOR
		\STATE Output the updated parameters $\boldsymbol{W}^k, \boldsymbol{W}^{T^k}$ and $\theta^k$ .
	\end{algorithmic}
\end{algorithm}

\subsection{Phase Retrieval}

The phase retrieval process involves recovering the phase information of the original signal  $\boldsymbol{x}$ from the measurements  $\boldsymbol{y}$
that only contain amplitude information, thereby obtaining the complete original signal. After adding the learning parameters, the model (\ref{model}) becomes
\begin{equation}
	\min_{\boldsymbol{x}} \frac{1}{2}\||\boldsymbol{W}\boldsymbol{x}|-\boldsymbol{y}\|_2^2+\lambda \|\boldsymbol{\varPsi}\boldsymbol{x}\|_1. \label{the phase retrieval with learned measurement matrix}
\end{equation}

The core idea for solving this problem is to combine the interpretability of the first-order proximal gradient optimization algorithm with the strong expressive power of neural networks. This is achieved by truncating the ISTA approximation and mapping it to a deep network structure with a fixed number of stages. Inspired by compressed sensing problems, we imitate the relevant steps of PRISTA-Net\cite{liu2023prista}; each stage consists of two modules: a subgradient descent module (SGD module) and a proximal mapping module (PPM module)

\textbf{Subgradient Descent Module (SGD Module)}

To solve Model (\ref{the phase retrieval with learned measurement matrix}), the first-order gradient descent method is usually adopted. Since the data fidelity term
$f(\boldsymbol{x}) = \frac{1}{2} \lVert |\boldsymbol{W} \boldsymbol{x}| - \boldsymbol{y} \rVert_2^2$ is non-differentiable, the subgradient of  $f(\boldsymbol{x})$
is used instead of the gradient, and its expression is
\begin{equation}
	\boldsymbol{W}^T \left( \boldsymbol{W} \boldsymbol{x} - \boldsymbol{y} \circ \frac{\boldsymbol{W} \boldsymbol{x}}{|\boldsymbol{W}\boldsymbol{x}|} \right) \in \partial_{\boldsymbol{x}} \frac{1}{2} \|  |\boldsymbol{W} \boldsymbol{x}| - \boldsymbol{y}\|_2^2,
\end{equation}

where the symbol $\circ$ denotes the element-wise Hadamard product, and 
$\partial_x \frac{1}{2}\|y - |\boldsymbol{W} \boldsymbol{x}|\|_2^2$ 
​is the subdifferential set of $f(\boldsymbol{W})$ at  $x$. Therefore, the SGD module can be expressed as
\begin{equation}
	\begin{split}
		\boldsymbol{r}^k &= \mathcal{SGD}(\boldsymbol{x}^{k-1}, t^k, \boldsymbol{y}, \boldsymbol{W}^{k-1}, {\boldsymbol{W}^T}^{k-1}) \\
		&= \boldsymbol{x}^{k-1} - t^k {\boldsymbol{W}^T}^{k-1}  \left( \boldsymbol{W}^{k-1} \boldsymbol{x}^{k-1} - \boldsymbol{y} \circ \frac{\boldsymbol{W}^{k-1} \boldsymbol{x}^{k-1}}{|\boldsymbol{W}^{k-1} \boldsymbol{x}^{k-1}|} \right).
	\end{split}
\end{equation}

\textbf{Proximal Projection Module} 
The PPM extracts sparse features through thresholding and non-linear transformation, thereby achieving signal sparsification and imposing sparse priors. The module is
\begin{equation}
	\boldsymbol{x}^k = \arg\min_{\boldsymbol{x}} \frac{1}{2} \| \boldsymbol{x} - \boldsymbol{r}^k\|_2^2 + \lambda \|\boldsymbol{\varPsi}\boldsymbol{x}\|_1.
\end{equation}

Elements of the intermediate variable  $\boldsymbol{r}^k$ 
are compared against a threshold: elements exceeding the threshold are shrunk, and irrelevant elements are set to zero, thereby achieving signal sparsification. Additionally, end-to-end learning is used to adaptively optimize the threshold, replacing the threshold in traditional methods that relies on cross-validation or empirical setting, which improves the quality of sparse solutions and convergence speed.

From Theorem 1 in reference \cite{zhang2018ista}, we have 
​ $\left\|\mathcal{F}(\boldsymbol{x}) - \mathcal{F}\left(\boldsymbol{r}^{(k)}\right) \right\|_2^2 \approx \alpha \left\| \boldsymbol{x} - \boldsymbol{r}^{k} \right\|_2^2$
. A parameterized learnable nonlinear transformation 
$\boldsymbol{\mathcal{F}}$ is used to replace $\boldsymbol{\varPsi}$
for imposing the sparse prior.
\begin{equation}
	\boldsymbol{x}^k = \arg\min_{\boldsymbol{x}} \frac{1}{2} \| \mathcal{F}(\boldsymbol{x}) - \mathcal{F}(\boldsymbol{r}^k)\|_2^2 + \eta^{k} \|\mathcal{F}(\boldsymbol{x})\|_1,
\end{equation}
here $\eta^k$ is a learnable threshold and let  $\boldsymbol{z}^{k} = \mathcal{F}(\boldsymbol{r}^k)$. Therefore, the closed-form solution of 
$\mathcal{F}(\boldsymbol{x}^k)$ 
can be derived as follows
\begin{equation}
	\tilde{\boldsymbol{z}}^{k} = \mathcal{F}(\boldsymbol{x}^k) = \text{soft}(\mathcal{F}(\boldsymbol{r}^k), \eta^k) = \text{soft}(\boldsymbol{z}^{k}, \eta^k).
\end{equation}

Inspired by the invertibility of transformation operators in traditional iterative algorithms, we denote the left inverse of  $\mathcal{F}$ as  $\tilde{\mathcal{F}}$, and thus we have
\begin{equation}
	\boldsymbol{x}^k = \tilde{\mathcal{F}} \left( \text{soft} \left( \mathcal{F}(\boldsymbol{r}^k), \eta^k \right) \right).
\end{equation}
To prevent gradient explosion or vanishing, the entire PPM module is designed as a residual structure
\begin{equation}
	\begin{split}
		\boldsymbol{x}^k &= \mathcal{P}\mathcal{P}\mathcal{M}_{\phi^k} \left(\boldsymbol{r}^k, \mathcal{F}^k, \tilde{\mathcal{F}}^k \right) \\
		&= \boldsymbol{r}^k + \tilde{\mathcal{F}}^k \left(\text{soft} \left(\mathcal{F}^{k} \left(\boldsymbol{r}^k \right),\eta^k\right)\right).
	\end{split}
\end{equation}
\begin{algorithm}[H]
	\caption{DLMMPR}
	\label{PRLISTA algorithm}
	\begin{algorithmic}[1]
		\REQUIRE Let $\boldsymbol{x}_0 $ be the initial point, and initialize the learnable parameters.
		\FOR{$ k=1,2,3,\ldots ,$K}
		\STATE $\boldsymbol{r}^k = \mathcal{SGD}(\boldsymbol{x}^{k-1}, t^k, \boldsymbol{y}, \boldsymbol{W}^{k-1}, {\boldsymbol{W}^T}^{k-1})$.
		\STATE $\boldsymbol{z}^{k}=\mathcal{F}\left(\boldsymbol{r}^{k}\right)$.
		\STATE $\tilde{\boldsymbol{z}}^{k} = \text{soft}\left( \boldsymbol{z}^k, \eta^k \right)$.
		\STATE$\boldsymbol{x}^{k}=\boldsymbol{r}^{k} + \tilde{\mathcal{F}}^k(\tilde{\boldsymbol{z}}^{k})$.
		\ENDFOR
		\STATE Output the current point $\boldsymbol{x}^{K}$ as the approximate optimal solution.
	\end{algorithmic}
\end{algorithm}

\begin{figure}[!h]
	\centering
	\renewcommand{\figurename}{Figure}
	\includegraphics[width=1\textwidth]{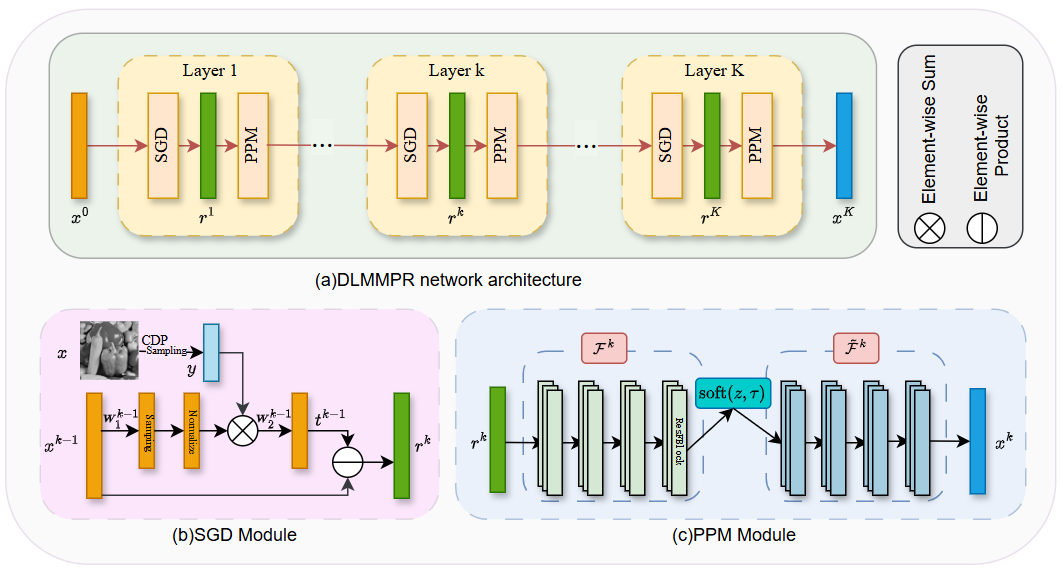}
	\caption{DLMMPR Network.}
	\label{fig:network} 
\end{figure}

We adopt the design from reference \cite{liu2023prista}, and design the nonlinear transformation function $\mathcal{F}$ into four parts. The network diagram of DLMMPR is shown in Figure \ref{fig:network}.

(1) 3×3 convolution: Perform 3×3 convolution operation on the intermediate result 
$\boldsymbol{r}^k$ to ensure high-throughput information transmission and avoid information loss.

(2) Convolutional Block Attention Module (CBAM): Adaptively learn the weights of different channels and spatial positions through channel attention and spatial attention modules, enabling the model to focus on key channel and spatial information. Through the organic combination of channel attention and spatial attention, the CBAM module can extract features with stronger representational ability and richer contextual information, thereby guiding the network to focus on phase information such as image edges, textures, and structures.

(3) Convolutional Unit (ConvU): This unit is composed of three 3×3 convolutional layers with ReLU activation functions, which is used to capture more detailed features.

(4) Fourier Residual Block (ResFBlock): This module enhances the ability to represent low-frequency and high-frequency features through the complementary fusion of global information and local details.

The optimized features are processed by a soft thresholding operator to filter out irrelevant information, thereby achieving signal sparsification. The left inverse transformation $\tilde{\mathcal{F}}$ adopts a symmetric structure design to 
$\mathcal{F}$, consisting of ResFBlock, ConvU, CBAM, and 3×3 convolution, which restores the corrected high-throughput information into high-quality reconstructed images. These designs have an enhancing effect on model performance: they improve the complementary learning of local and global information, enhance the focus on high-dimensional features including image edges, textures, and structures, and the designed logarithmic loss achieves significant improvement at low noise levels.

\subsection{Parameters and Initialization}
The set of learnable parameters for the proposed architecture is denoted as 
$\Theta = \{\boldsymbol{W}^k, {\boldsymbol{W}^T}^{k}, t^k$,
$ \eta^k, \mathcal{F}^k, \tilde{\mathcal{F}}^k\}_{k=1}^{K}$, which includes the measurement matrix 
$\boldsymbol{W}^k_1$ in the SGD module and its transpose matrix 
$\boldsymbol{W}^k_2$, the iteration step size 
$t^k$, the threshold parameter $\eta^k$
in the PPM module, as well as the nonlinear transformation function 
$\mathcal{F}^k$ and its corresponding inverse transformation 
$\tilde{\mathcal{F}}^k$. All parameters can be learned end-to-end by minimizing the loss function. The network initialization adopts the Xavier initialization method \cite{2010Understanding}, where the initial values of 
$\eta^1 $ and $\theta^1$ 
are set to 0.5 and 0.01, respectively. The default configuration of the architecture includes 7 processing stages, with the number of convolution channels in each stage set to 32, and parameters are not shared between different stages. Following the initialization strategy of PrDeep, the initial iteration input is set as 
$x^0 = 1$.
\section{Numerical Experiments}
The algorithm was implemented in PyTorch, employing the Adam optimizer with an initial learning rate of $1 \times 10^{-3}$. A step-wise learning rate decay of $5\%$
was applied every two epochs, and a batch size of 10 was used. Experiments were performed on a workstation comprising an Intel(R) Core(TM) Ultra 9 285K @3.70 GHz CPU and an NVIDIA RTX A5000 (24GB) GPU. Investigations into the number of network stages revealed suboptimal performance when set between 3 and 6 stages. Beyond 7 stages, performance gains became negligible while computational runtime escalated considerably. Accordingly,  we selected 7 as the optimal number of network stages.
and trained for 150 epochs, a process requiring approximately 5 hours.
\subsection{Experimental Setup}
\begin{figure}[t]
	\centering
	\captionsetup[subfigure]{skip=0.1pt} 
	\hspace{-\tabcolsep} 
	
	\begin{subfigure}{0.15\linewidth}
		\centering
		\includegraphics[width=\linewidth]{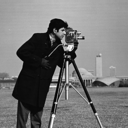}
		\footnotesize{Cameraman}
	\end{subfigure}%
	\hspace{1pt}
	\begin{subfigure}{0.15\linewidth}
		\centering
		\includegraphics[width=\linewidth]{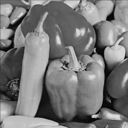}
		\footnotesize{Peppers}
	\end{subfigure}%
	\hspace{1pt}
	\begin{subfigure}{0.15\linewidth}
		\centering
		\includegraphics[width=\linewidth]{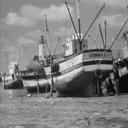}
		\footnotesize{Boat}
	\end{subfigure}%
	\hspace{1pt}
	\begin{subfigure}{0.15\linewidth}
		\centering
		\includegraphics[width=\linewidth]{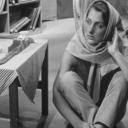}
		\footnotesize{Barbara}
	\end{subfigure}%
	\hspace{1pt}
	\begin{subfigure}{0.15\linewidth}
		\centering
		\includegraphics[width=\linewidth]{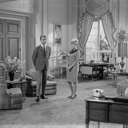}
		\footnotesize{Couple}
	\end{subfigure}%
	\hspace{1pt}
	\begin{subfigure}{0.15\linewidth}
		\centering
		\includegraphics[width=\linewidth]{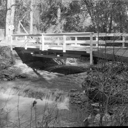}
		\footnotesize{Streamandbridge}
	\end{subfigure}
	
	\vspace{6pt}
	
	\begin{subfigure}{0.15\linewidth}
		\centering
		\includegraphics[width=\linewidth]{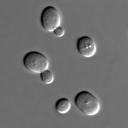}
		\footnotesize{Yeast}
	\end{subfigure}%
	\hspace{1pt}
	\begin{subfigure}{0.15\linewidth}
		\centering
		\includegraphics[width=\linewidth]{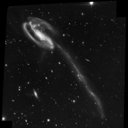}
		\footnotesize{Tadpole Galaxy}
	\end{subfigure}%
	\hspace{1pt}
	\begin{subfigure}{0.15\linewidth}
		\centering
		\includegraphics[width=\linewidth]{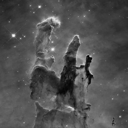}
		\footnotesize{Pillars Creation}
	\end{subfigure}%
	\hspace{1pt}
	\begin{subfigure}{0.15\linewidth}
		\centering
		\includegraphics[width=\linewidth]{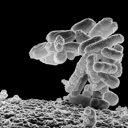}
		\footnotesize{Ecoli}
	\end{subfigure}%
	\hspace{1pt}
	\begin{subfigure}{0.15\linewidth}
		\centering
		\includegraphics[width=\linewidth]{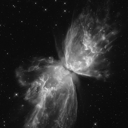}
		\footnotesize{Butterfly}
	\end{subfigure}%
	\hspace{1pt}
	\begin{subfigure}{0.15\linewidth}
		\centering
		\includegraphics[width=\linewidth]{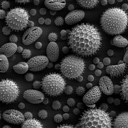}
		\footnotesize{Pollen}
	\end{subfigure}
	\caption{Widely Used Test Datasets in Phase Retrieval Set12.}
	\label{Set12}
\end{figure}

For training purposes, this study leveraged a comprehensive dataset comprising 6000 grayscale images, synthesized from the publicly available BSD400 \cite{BSD400} and PASCAL VOC \cite{PASCAL} datasets. Each image underwent uniform resizing to a consistent resolution of 128×128 pixels. These prepared images served as the foundation for generating the measurement values (CDPs), following the procedure outlined in Equation (\ref{measurement}). To rigorously evaluate our model’s performance, we adopted the test dataset employed by PrDeep \cite{metzler2018prdeep}. This comprises two curated sets: the NT-6 set, featuring 6 natural images, and the UNT-6 set, containing 6 unnatural images, both widely recognized in phase retrieval literature (detailed in Figure \ref{Set12}). Our experimental focus remained predominantly on the 
resolution scale, and all subsequent method comparisons were strictly performed at this uniform resolution.

For seamless integration with model input requirements, all image pixel values were uniformly normalized to the range [0, 1]. The intensity of Poisson noise was calibrated to the pixel range by dividing the noise level by 255. During the training regimen, the noise level parameter, denoted as  $w$
, was drawn uniformly from the discrete set $\{9, 27, 81\}$. To rigorously assess the algorithm’s practical generalization capabilities, the test  employed a measurement mask distribution,  $D$ , that differed from that utilized during training.
\subsection{Model Evaluation}
1. \textbf{Quantitative Evaluation}   
~Across noise levels $\alpha=9$, $27$ and $81$, the model consistently achieves optimal or near-optimal reconstruction performance by 150 epochs, as shown in Table \ref{tab:epoch_performance}
. In the low-noise case, UNT-6 PSNR at 150 epochs ($41.887$dB) marginally exceeds later peak performance ($41.864$dB). Similar trends hold for noise levels 27 and 81. The performance saturates around 150 epochs, with subsequent training yielding little benefit.  Moreover, the higher the noise level, the more obvious this effect,  showing the characteristics of early convergence and high stability.

\begin{table}[t]
	\centering
	\renewcommand{\arraystretch}{0.8}
	\caption{Performance Performance Across Different Noise Levels and Training }
	\label{tab:epoch_performance}
	\begin{tabular}{ccccccc} 
		\toprule 
		\multirow{2}{*}{epoch} & \multicolumn{2}{c}{$\alpha = 9$} & \multicolumn{2}{c}{$\alpha = 27$} & \multicolumn{2}{c}{$\alpha = 81$} \\
		\cmidrule(lr){2-3} \cmidrule(lr){4-5} \cmidrule(lr){6-7} 
		& UNT-6 & NT-6 & UNT-6 & NT-6 & UNT-6 & NT-6 \\
		\midrule[\heavyrulewidth] 
		25 & 41.76/0.98 & 40.27/0.99 & 35.09/0.94 & 33.39/0.94 & 29.44/0.85 & 27.30/0.83 \\
		50 & 41.82/0.98 & 40.35/0.98 & 35.14/0.94 & 33.41/0.95 & 29.52/0.84 & 27.30/0.83 \\
		75 & 41.83/0.98 & 40.38/0.99 & 35.18/0.95 & 33.47/0.95 & 29.61/0.86 & 27.36/0.84 \\
		100 & 41.85/0.98 & 40.37/0.99 & 35.20/0.95 & 33.50/0.95 & 29.65/0.86 & 27.36/0.83 \\
		125 & 41.85/0.98 & 40.40/0.99& 35.21/0.94 & 33.48/0.95 & \textbf{29.74/0.86} & 27.40/0.83 \\
		150 & \textbf{41.89/0.98} & 40.41/0.99 & \textbf{35.22/0.95} & \underline{33.51}/0.95 & 29.70/0.86 & \textbf{27.44/0.83} \\
		175 & 41.86/0.98 & \underline{40.41}/0.99 & 35.19/0.95 & 33.48/0.95 & 29.68/0.86 & 27.38/0.83 \\
		200 & 41.85/0.98 & 40.39/0.99 & \underline{35.22}/0.95 & \textbf{33.52/0.95} & \underline{29.70}/0.86 & \underline{27.43}/0.85 \\
		225 & \underline{41.86}/0.98 & 40.37/0.99 & 35.19/0.95 & 33.48/0.95 & 29.64/0.86 & 27.33/0.83 \\
		250 & 41.86/0.98 & \textbf{40.42/0.99} & 35.17/0.95 & 33.41/0.95 & 29.61/0.86 & 27.36/0.83 \\
		275 & 41.82/0.98 & 40.38/0.99 & 35.16/0.95 & 33.47/0.95 & 29.63/0.86 & 27.36/0.83 \\
		300 & 41.86/0.98 & 40.39/0.99 & 35.17/0.95 & 33.50/0.95 & 29.63/0.86 & 27.38/0.83 \\
		\bottomrule 
	\end{tabular}
\end{table}

The presented results compellingly underscore the power of learning to optimize measurement matrix parameters. Optimizing matrix  $A$ through learning, rather than relying on fixed or hand-designed matrices, demonstrably reduces training epochs and accelerates convergence. DLMMPR thus reaches its optimal solution earlier in iterations, avoiding the slow, performance-limited convergence of traditional methods. The specific numerical differences intuitively reflect the profound improvement learned parameters bring to phase retrieval.

\begin{figure}[htbp]
	\centering
	\captionsetup[subfigure]{skip=0.1pt} 
	\hspace{-\tabcolsep} 
	
	\begin{subfigure}{0.15\linewidth}
		\centering
		\includegraphics[width=\linewidth]{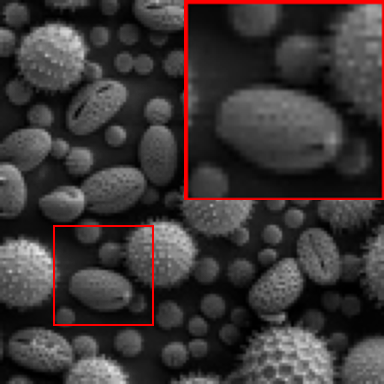}
		\footnotesize{25/38.691}
	\end{subfigure}%
	\hspace{1pt}
	\begin{subfigure}{0.15\linewidth}
		\centering
		\includegraphics[width=\linewidth]{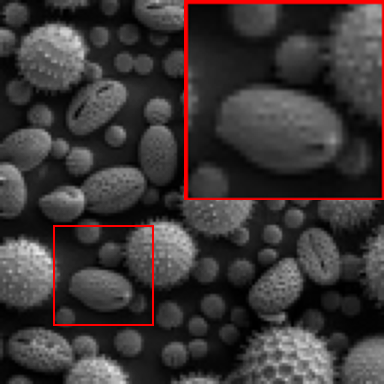}
		\footnotesize{50/38.742}
	\end{subfigure}%
	\hspace{1pt}
	\begin{subfigure}{0.15\linewidth}
		\centering
		\includegraphics[width=\linewidth]{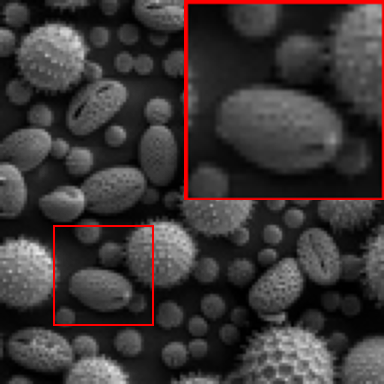}
		\footnotesize{75/38.819}
	\end{subfigure}%
	\hspace{1pt}
	\begin{subfigure}{0.15\linewidth}
		\centering
		\includegraphics[width=\linewidth]{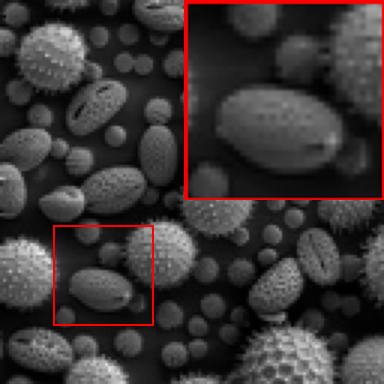}
		\footnotesize{100/38.795}
	\end{subfigure}%
	\hspace{1pt}
	\begin{subfigure}{0.15\linewidth}
		\centering
		\includegraphics[width=\linewidth]{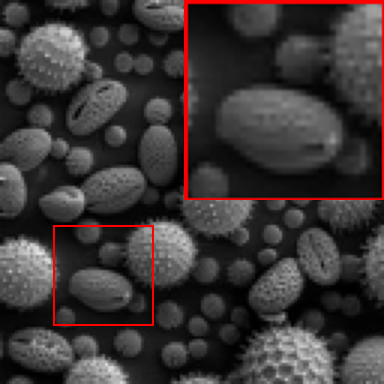}
		\footnotesize{125/38.757}
	\end{subfigure}%
	\hspace{1pt}
	\begin{subfigure}{0.15\linewidth}
		\centering
		\includegraphics[width=\linewidth]{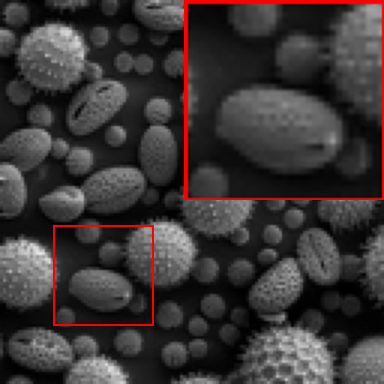}
		\footnotesize{150/38.877}
	\end{subfigure}
	
	\vspace{6pt}
	
	\begin{subfigure}{0.15\linewidth}
		\centering
		\includegraphics[width=\linewidth]{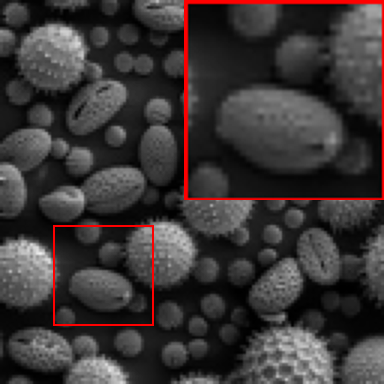}
		\footnotesize{175/38.861}
	\end{subfigure}%
	\hspace{1pt}
	\begin{subfigure}{0.15\linewidth}
		\centering
		\includegraphics[width=\linewidth]{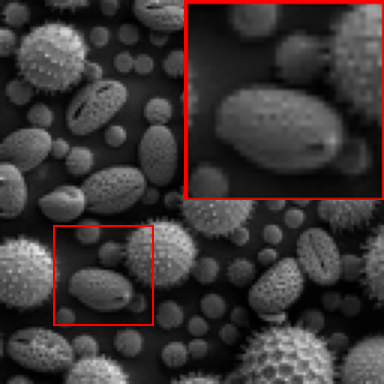}
		\footnotesize{200/38.775}
	\end{subfigure}%
	\hspace{1pt}
	\begin{subfigure}{0.15\linewidth}
		\centering
		\includegraphics[width=\linewidth]{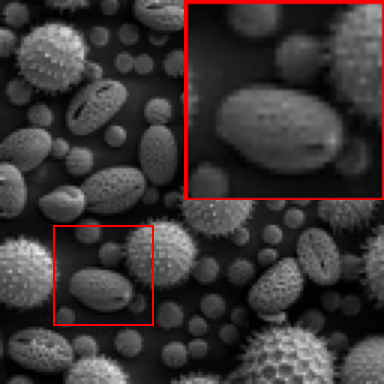}
		\footnotesize{225/38.790}
	\end{subfigure}%
	\hspace{1pt}
	\begin{subfigure}{0.15\linewidth}
		\centering
		\includegraphics[width=\linewidth]{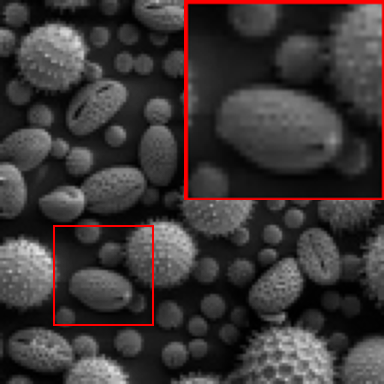}
		\footnotesize{250Pollen/38.817}
	\end{subfigure}%
	\hspace{1pt}
	\begin{subfigure}{0.15\linewidth}
		\centering
		\includegraphics[width=\linewidth]{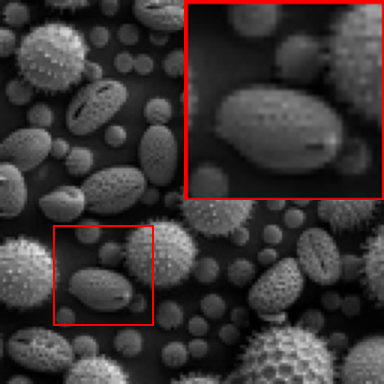}
		\footnotesize{275/38.772}
	\end{subfigure}%
	\hspace{1pt}
	\begin{subfigure}{0.15\linewidth}
		\centering
		\includegraphics[width=\linewidth]{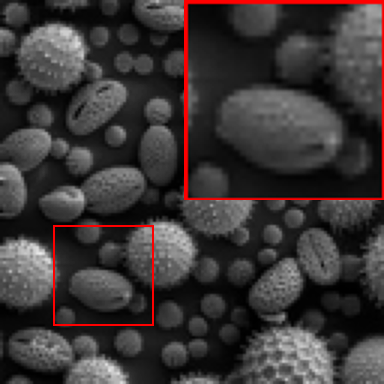}
		\footnotesize{300/38.867}
	\end{subfigure}
	\caption{Phase Retrieval Results of Pollen Image Under Different Iteration Counts at Noise Level  $\alpha = 9$ epoch/PSNR(dB).}
	\label{fig:pollen_phase_retrieval}
\end{figure}

2. \textbf{Qualitative Evaluation}  
~Figure \ref{fig:pollen_phase_retrieval}
illustrates the phase retrieval progression for the Pollen image across discrete 25-iteration intervals, up to a total of 300 iterations, under a noise level of  $\alpha = 9$. A thorough examination of both the visual outcomes and quantitative metrics unequivocally demonstrates that the proposed method’s reconstruction performance exhibits a dynamic optimization trajectory throughout the iterative process, culminating in an optimal state at approximately 150 iterations.

Visually, the 150-iteration reconstruction excels by balancing texture detail and noise suppression. Edge contours of pollen grains are sharp, and surface texture layering is superior to other iterations, with fine grooves and convex structures fully restored. Earlier iterations (25–75) show blurriness and halos (at 75 iterations). While 100–125 iterations improve detail, their PSNR lags by 0.08–0.12 dB. Later iterations (175–300) offer comparable PSNR but exhibit subtle over-smoothing, reducing clarity of fine grain surface lines.

The decision to set iteration count at 150 is therefore justified. This ensures complete learning of image structural features and texture details while preventing accuracy loss from excessive iteration. It is the optimal choice under current experimental conditions, confirming DLMMPR’s controllability and stable reconstruction quality during iteration.
\subsection{Comparison with State-of-the-Art Methods}
In advancing the field of phase retrieval algorithms, we benchmark our proposed method against leading contemporary techniques. Our selection for comparison includes DeepMMSE \cite{chen2022unsupervised}, PrComplex \cite{chen2022phase}, PRISTA-Net \cite{liu2023prista}, and USP \cite{quan2023unsupervised}. To establish a robust and equitable comparison, we leveraged the publicly disseminated code for each algorithm and conducted experiments on the same public test dataset used for our evaluation. This approach ensures that performance differences are attributable to the algorithmic methodologies themselves, rather than variations in experimental conditions.


\begin{table}[htbp]
	\centering
	\renewcommand{\arraystretch}{0.8}
	\small
	\caption{Phase Retrieval Performance Comparison on the prDeep Set12 Dataset }
	\label{tab:pr_comparison}
	\begin{tabular}{@{}l@{\hspace{1.3\tabcolsep}}c@{\hspace{1.3\tabcolsep}}c@{\hspace{1.3\tabcolsep}}c@{\hspace{1.3\tabcolsep}}c@{\hspace{1.3\tabcolsep}}c@{\hspace{1.3\tabcolsep}}c@{}}
		\toprule
		\multirow{2}{*}{\vspace{-0.5em}Method} & \multicolumn{2}{c}{$\alpha = 9$} & \multicolumn{2}{c}{$\alpha = 27$} & \multicolumn{2}{c}{$\alpha = 81$} \\
		\cmidrule(lr){2-3} \cmidrule(lr){4-5} \cmidrule(lr){6-7}
		& UNT-6 & NT-6 & UNT-6 & NT-6 & UNT-6 & NT-6 \\
		\midrule
		DeepMMSE\cite{chen2022unsupervised} & 40.26/0.98 & 39.45/0.98 & 33.03/0.94 & 32.32/0.93 & 26.64/0.81 & 25.41/0.78 \\
		prCom\cite{chen2022phase} & 41.11/0.98 & 39.78/0.99 & 34.64/0.93 & 33.49/0.94 & 28.27/0.84 & 26.55/0.82 \\
		PRISTA-Net\cite{liu2023prista} & 41.47/0.98 & 40.10/0.99 & 34.89/0.94 & 33.31/0.94 & 29.45/0.85 & 27.08/0.82 \\
		UPR\cite{quan2023unsupervised} & \underline{41.75/0.98} & \underline{40.26/0.99} & \underline{34.98/0.94} & \underline{33.33/0.94} & \underline{29.57/0.86} & \underline{27.19/0.83} \\
		DLMMPR & \textbf{41.89/0.98} & \textbf{40.41/0.99} & \textbf{35.23/0.95} & \textbf{33.51/0.95} & \textbf{29.70/0.86} & \textbf{27.44/0.83} \\
		\bottomrule
	\end{tabular}
\end{table}

\begin{table}[h!]
	\centering
	\renewcommand{\arraystretch}{0.8} 
	\setlength{\tabcolsep}{7pt} 
	\caption{Performance Comparison of Different Images Under Various Noise Levels }
	\label{tab:dataset_performance}
	\begin{tabular}{ccccccc} 
		\toprule 
		& $\alpha$ & DMMSE & PrCom & PRISTA-Net & USP & Ours \\
		\midrule 
		\multirow{3}{*}{Ecoli} & 9 & 39.91/0.996 & 39.80/0.967 & 39.77/0.989 & \underline{40.74/0.994} & \textbf{40.80/0.995} \\
		& 27 & 31.02/0.971 & 32.35/0.914 & 32.36/0.964 & \underline{32.44/0.967} & \textbf{32.49/0.968} \\
		& 81 & 23.49/0.850 & 25.21/0.786 & 25.55/0.866 & \underline{25.28/0.884} & \textbf{25.68/0.902} \\
		\midrule 
		\multirow{3}{*}{Yeast} & 9 & 42.90/0.987 & 43.95/0.983 & 45.54/0.990 & \underline{45.67/0.990} & \textbf{45.79/0.990} \\
		& 27 & 37.06/0.975 & 38.36/0.963 & 39.71/0.981 & \underline{39.95/0.983} & \textbf{40.15/0.982} \\
		& 81 & 27.52/0.900 & 29.41/0.922 & 33.12/0.954 & \underline{33.63/0.960} & \textbf{33.75/0.958} \\
		\midrule 
		\multirow{3}{*}{Butterfly} & 9 & 40.47/0.97 & \underline{41.63/0.977} & 41.53/0.998 & 41.63/0.978 & \textbf{41.81/0.979} \\
		& 27 & 33.73/0.91 & 35.17/0.928 & 35.18/0.929 & \underline{35.31/0.933} & \textbf{35.55/0.934} \\
		& 81 & 27.70/0.792 & 29.10/0.843 & \underline{30.38/0.867} & 30.20/0.841 & \textbf{30.68/0.853} \\
		\midrule 
		\multirow{3}{*}{barbara} & 9 & 39.82/0.986 & 40.03/0.987 & 40.66/0.989 & \underline{40.83/0.990} & \textbf{40.96/0.990} \\
		& 27 & 32.77/0.95 & \underline{34.12/0.960} & 34.05/0.961 & 34.03/0.961 & \textbf{34.43/0.964} \\
		& 81 & 25.52/0.833 & 26.67/0.868 & 27.69/0.875 & \underline{27.73/0.881} & \textbf{28.05/0.878} \\
		\midrule 
		\multirow{3}{*}{cameraman} & 9 & 40.18/0.975 & 41.08/0.980 & 41.34/0.983 & \underline{41.37/0.983} & \textbf{41.60/0.984} \\
		& 27 & 33.24/0.917 & 34.66/0.942 & 34.857/0.951 & \underline{34.912/0.951} & \textbf{35.04/0.954} \\
		& 81 & 26.10/0.815 & 27.75/0.875 & 28.36/0.875 & \underline{28.75/0.874} & \textbf{28.95/0.878} \\
		\midrule 
		\multirow{3}{*}{peppers}  & 9  & 40.22/0.985 & 40.20/0.985 & 41.07/0.988 & \underline{41.14/0.988} & \textbf{41.27/0.989} \\
		& 27 & 33.27/0.949 & 34.23/0.960 & \underline{34.44/0.961} & 34.40/0.960 & \textbf{34.55/0.961} \\
		& 81 & 25.61/0.834 & 26.79/0.871 & 27.80/0.870 & \underline{27.99/0.878} & \textbf{28.10/0.880} \\
		\bottomrule 
	\end{tabular}
\end{table}
\begin{figure}[!t]
	\centering
	\captionsetup[subfigure]{skip=0.1pt} 
	\hspace{-\tabcolsep} 
	
	\begin{subfigure}{0.15\linewidth}
		\centering
		{Truth}
		\includegraphics[width=\linewidth]{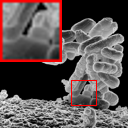}
		\small{PSNR/SSIM}
	\end{subfigure}%
	\hspace{1pt}
	\begin{subfigure}{0.15\linewidth}
		\centering
		{DeepMMSE}
		\includegraphics[width=\linewidth]{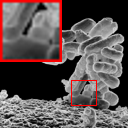}
		\small{39.91/0.996}
	\end{subfigure}%
	\hspace{1pt}
	\begin{subfigure}{0.15\linewidth}
		\centering
		{PrCom}
		\includegraphics[width=\linewidth]{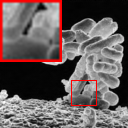}
		\small{39.80/0.982}
	\end{subfigure}%
	\hspace{1pt}
	\begin{subfigure}{0.15\linewidth}
		\centering
		{PRISTA-Net}
		\includegraphics[width=\linewidth]{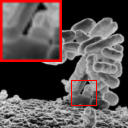}
		\small{39.77/0.988}
	\end{subfigure}%
	\hspace{1pt}
	\begin{subfigure}{0.15\linewidth}
		\centering
		{USP}
		\includegraphics[width=\linewidth]{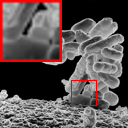}
		\small{40.74/0.993}
	\end{subfigure}%
	\hspace{1pt}
	\begin{subfigure}{0.15\linewidth}
		\centering
		{Ours}
		\includegraphics[width=\linewidth]{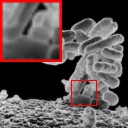}
		\small{40.80/0.995}
	\end{subfigure}
	
	\vspace{6pt}
	\begin{subfigure}{0.15\linewidth}
		\centering
		\includegraphics[width=\linewidth]{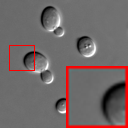}
		\footnotesize{PSNR/SSIM}
	\end{subfigure}%
	\hspace{1pt}
	\begin{subfigure}{0.15\linewidth}
		\centering
		\includegraphics[width=\linewidth]{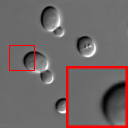}
		\small{42.90/0.987}
	\end{subfigure}%
	\hspace{1pt}
	\begin{subfigure}{0.15\linewidth}
		\centering
		\includegraphics[width=\linewidth]{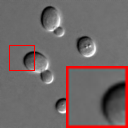}
		\small{43.95/0.982}
	\end{subfigure}%
	\hspace{1pt}
	\begin{subfigure}{0.15\linewidth}
		\centering
		\includegraphics[width=\linewidth]{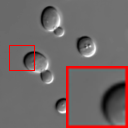}
		\small{45.54/0.989}
	\end{subfigure}%
	\hspace{1pt}
	\begin{subfigure}{0.15\linewidth}
		\centering
		\includegraphics[width=\linewidth]{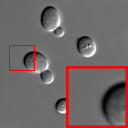}
		\small{45.67/0.989}
	\end{subfigure}%
	\hspace{1pt}
	\begin{subfigure}{0.15\linewidth}
		\centering
		\includegraphics[width=\linewidth]{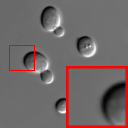}
		\small{45.79/0.990}
	\end{subfigure}
	
	\vspace{6pt}
	
	\begin{subfigure}{0.15\linewidth}
		\centering
		\includegraphics[width=\linewidth]{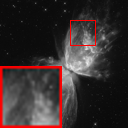}
		\small{PSNR/SSIM}
	\end{subfigure}%
	\hspace{1pt}
	\begin{subfigure}{0.15\linewidth}
		\centering
		\includegraphics[width=\linewidth]{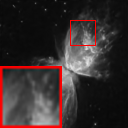}
		\small{40.47/0.970}
	\end{subfigure}%
	\hspace{1pt}
	\begin{subfigure}{0.15\linewidth}
		\centering
		\includegraphics[width=\linewidth]{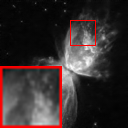}
		\small{41.63/0.976}
	\end{subfigure}%
	\hspace{1pt}
	\begin{subfigure}{0.15\linewidth}
		\centering
		\includegraphics[width=\linewidth]{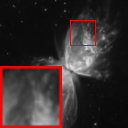}
		\small{41.53/0.977}
	\end{subfigure}%
	\hspace{1pt}
	\begin{subfigure}{0.15\linewidth}
		\centering
		\includegraphics[width=\linewidth]{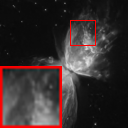}
		\small{41.63/0.978}
	\end{subfigure}%
	\hspace{1pt}
	\begin{subfigure}{0.15\linewidth}
		\centering
		\includegraphics[width=\linewidth]{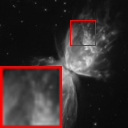}
		\small{41.81/0.979}
	\end{subfigure}
	
	\vspace{6pt}
	\begin{subfigure}{0.15\linewidth}
		\centering
		\includegraphics[width=\linewidth]{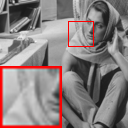}
		\small{PSNR/SSIM}
	\end{subfigure}%
	\hspace{1pt}
	\begin{subfigure}{0.15\linewidth}
		\centering
		\includegraphics[width=\linewidth]{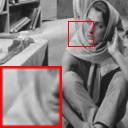}
		\small{39.82/0.986}
	\end{subfigure}%
	\hspace{1pt}
	\begin{subfigure}{0.15\linewidth}
		\centering
		\includegraphics[width=\linewidth]{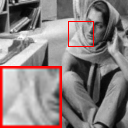}
		\small{40.03/0.987}
	\end{subfigure}%
	\hspace{1pt}
	\begin{subfigure}{0.15\linewidth}
		\centering
		\includegraphics[width=\linewidth]{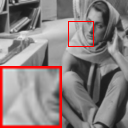}
		\small{40.66/0.989}
	\end{subfigure}%
	\hspace{1pt}
	\begin{subfigure}{0.15\linewidth}
		\centering
		\includegraphics[width=\linewidth]{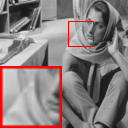}
		\small{40.83/0.990}
	\end{subfigure}%
	\hspace{1pt}
	\begin{subfigure}{0.15\linewidth}
		\centering
		\includegraphics[width=\linewidth]{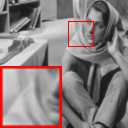}
		\small{40.96/0.990}
	\end{subfigure}
	
	\vspace{6pt}
	\begin{subfigure}{0.15\linewidth}
		\centering
		\includegraphics[width=\linewidth]{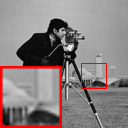}
		\small{PSNR/SSIM}
	\end{subfigure}%
	\hspace{1pt}
	\begin{subfigure}{0.15\linewidth}
		\centering
		\includegraphics[width=\linewidth]{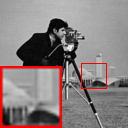}
		\small{40.18/0.975}
	\end{subfigure}%
	\hspace{1pt}
	\begin{subfigure}{0.15\linewidth}
		\centering
		\includegraphics[width=\linewidth]{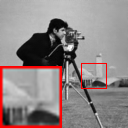}
		\small{41.08/0.980}
	\end{subfigure}%
	\hspace{1pt}
	\begin{subfigure}{0.15\linewidth}
		\centering
		\includegraphics[width=\linewidth]{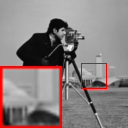}
		\small{41.34/0.983}
	\end{subfigure}%
	\hspace{1pt}
	\begin{subfigure}{0.15\linewidth}
		\centering
		\includegraphics[width=\linewidth]{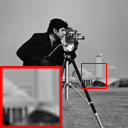}
		\small{41.37/0.983}
	\end{subfigure}%
	\hspace{1pt}
	\begin{subfigure}{0.15\linewidth}
		\centering
		\includegraphics[width=\linewidth]{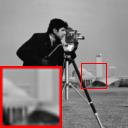}
		\small{41.60/0.984}
	\end{subfigure}
	
	\vspace{6pt}
	\begin{subfigure}{0.15\linewidth}
		\centering
		\includegraphics[width=\linewidth]{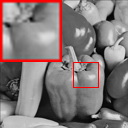}
		\small{PSNR/SSIM}
	\end{subfigure}%
	\hspace{1pt}
	\begin{subfigure}{0.15\linewidth}
		\centering
		\includegraphics[width=\linewidth]{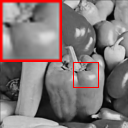}
		\small{40.22/0.985}
	\end{subfigure}%
	\hspace{1pt}
	\begin{subfigure}{0.15\linewidth}
		\centering
		\includegraphics[width=\linewidth]{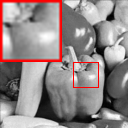}
		\small{40.20/0.985}
	\end{subfigure}%
	\hspace{1pt}
	\begin{subfigure}{0.15\linewidth}
		\centering
		\includegraphics[width=\linewidth]{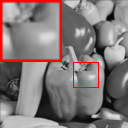}
		\small{41.07/0.988}
	\end{subfigure}%
	\hspace{1pt}
	\begin{subfigure}{0.15\linewidth}
		\centering
		\includegraphics[width=\linewidth]{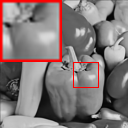}
		\small{41.14/0.988}
	\end{subfigure}%
	\hspace{1pt}
	\begin{subfigure}{0.15\linewidth}
		\centering
		\includegraphics[width=\linewidth]{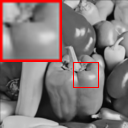}
		\small{41.27/0.989}
	\end{subfigure}
	
	\caption{Comparison of Image Restoration Effects of Different Methods at Noise Level $\alpha =9$.}
	\label{fig:phase_retrieval_comparison9}
\end{figure}

\begin{figure}[!t]
	\centering
	\captionsetup[subfigure]{skip=0.1pt} 
	\hspace{-\tabcolsep} 
	
	\begin{subfigure}{0.15\linewidth}
		\centering
		{Truth}
		\includegraphics[width=\linewidth]{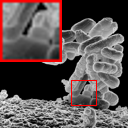}
		\small{PSNR/SSIM}
	\end{subfigure}%
	\hspace{1pt}
	\begin{subfigure}{0.15\linewidth}
		\centering
		{DeepMMSE}
		\includegraphics[width=\linewidth]{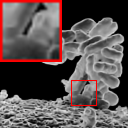}
		\small{31.02/0.971}
	\end{subfigure}%
	\hspace{1pt}
	\begin{subfigure}{0.15\linewidth}
		\centering
		{PrCom}
		\includegraphics[width=\linewidth]{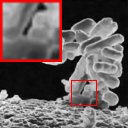}
		\small{32.35/0.914}
	\end{subfigure}%
	\hspace{1pt}
	\begin{subfigure}{0.15\linewidth}
		\centering
		{PRISTA-Net}
		\includegraphics[width=\linewidth]{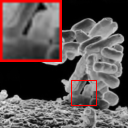}
		\small{32.36/0.964}
	\end{subfigure}%
	\hspace{1pt}
	\begin{subfigure}{0.15\linewidth}
		\centering
		{USP}
		\includegraphics[width=\linewidth]{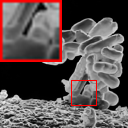}
		\small{32.44/0.967}
	\end{subfigure}%
	\hspace{1pt}
	\begin{subfigure}{0.15\linewidth}
		\centering
		{Ours}
		\includegraphics[width=\linewidth]{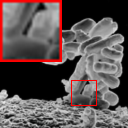}
		\small{32.49/0.968}
	\end{subfigure}
	
	\vspace{6pt}
	
	\begin{subfigure}{0.15\linewidth}
		\centering
		\includegraphics[width=\linewidth]{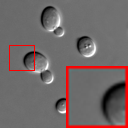}
		\small{PSNR/SSIM}
	\end{subfigure}%
	\hspace{1pt}
	\begin{subfigure}{0.15\linewidth}
		\centering
		\includegraphics[width=\linewidth]{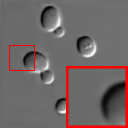}
		\small{37.06/0.975}
	\end{subfigure}%
	\hspace{1pt}
	\begin{subfigure}{0.15\linewidth}
		\centering
		\includegraphics[width=\linewidth]{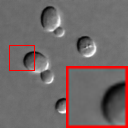}
		\small{38.36/0.963}
	\end{subfigure}%
	\hspace{1pt}
	\begin{subfigure}{0.15\linewidth}
		\centering
		\includegraphics[width=\linewidth]{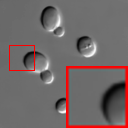}
		\small{39.71/0.981}
	\end{subfigure}%
	\hspace{1pt}
	\begin{subfigure}{0.15\linewidth}
		\centering
		\includegraphics[width=\linewidth]{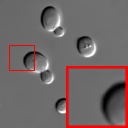}
		\small{39.95/0.983}
	\end{subfigure}%
	\hspace{1pt}
	\begin{subfigure}{0.15\linewidth}
		\centering
		\includegraphics[width=\linewidth]{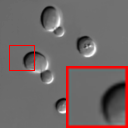}
		\small{40.15/0.982}
	\end{subfigure}
	
	\vspace{6pt}
	
	\begin{subfigure}{0.15\linewidth}
		\centering
		\includegraphics[width=\linewidth]{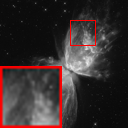}
		\small{PSNR/SSIM}
	\end{subfigure}%
	\hspace{1pt}
	\begin{subfigure}{0.15\linewidth}
		\centering
		\includegraphics[width=\linewidth]{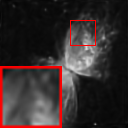}
		\small{33.73/0.910}
	\end{subfigure}%
	\hspace{1pt}
	\begin{subfigure}{0.15\linewidth}
		\centering
		\includegraphics[width=\linewidth]{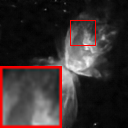}
		\small{35.17/0.928}
	\end{subfigure}%
	\hspace{1pt}
	\begin{subfigure}{0.15\linewidth}
		\centering
		\includegraphics[width=\linewidth]{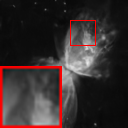}
		\small{35.18/0.929}
	\end{subfigure}%
	\hspace{1pt}
	\begin{subfigure}{0.15\linewidth}
		\centering
		\includegraphics[width=\linewidth]{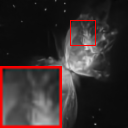}
		\small{35.31/0.933}
	\end{subfigure}%
	\hspace{1pt}
	\begin{subfigure}{0.15\linewidth}
		\centering
		\includegraphics[width=\linewidth]{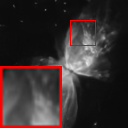}
		\small{35.55/0.934}
	\end{subfigure}
	
	\vspace{6pt}
	
	\begin{subfigure}{0.15\linewidth}
		\centering
		\includegraphics[width=\linewidth]{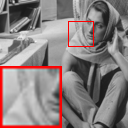}
		\small{PSNR/SSIM}
	\end{subfigure}%
	\hspace{1pt}
	\begin{subfigure}{0.15\linewidth}
		\centering
		\includegraphics[width=\linewidth]{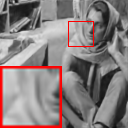}
		\small{32.77/0.950}
	\end{subfigure}%
	\hspace{1pt}
	\begin{subfigure}{0.15\linewidth}
		\centering
		\includegraphics[width=\linewidth]{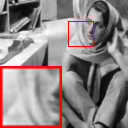}
		\small{34.12/0.960}
	\end{subfigure}%
	\hspace{1pt}
	\begin{subfigure}{0.15\linewidth}
		\centering
		\includegraphics[width=\linewidth]{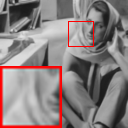}
		\small{34.05/0.961}
	\end{subfigure}%
	\hspace{1pt}
	\begin{subfigure}{0.15\linewidth}
		\centering
		\includegraphics[width=\linewidth]{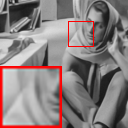}
		\small{34.03/0.961}
	\end{subfigure}%
	\hspace{1pt}
	\begin{subfigure}{0.15\linewidth}
		\centering
		\includegraphics[width=\linewidth]{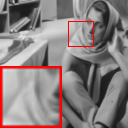}
		\small{34.43/0.964}
	\end{subfigure}
	
	\vspace{6pt}
	
	\begin{subfigure}{0.15\linewidth}
		\centering
		\includegraphics[width=\linewidth]{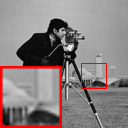}
		\small{PSNR/SSIM}
	\end{subfigure}%
	\hspace{1pt}
	\begin{subfigure}{0.15\linewidth}
		\centering
		\includegraphics[width=\linewidth]{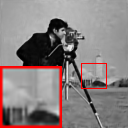}
		\small{33.24/0.917}
	\end{subfigure}%
	\hspace{1pt}
	\begin{subfigure}{0.15\linewidth}
		\centering
		\includegraphics[width=\linewidth]{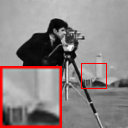}
		\small{34.66/0.942}
	\end{subfigure}%
	\hspace{1pt}
	\begin{subfigure}{0.15\linewidth}
		\centering
		\includegraphics[width=\linewidth]{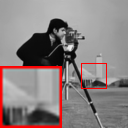}
		\small{34.85/0.951}
	\end{subfigure}%
	\hspace{1pt}
	\begin{subfigure}{0.15\linewidth}
		\centering
		\includegraphics[width=\linewidth]{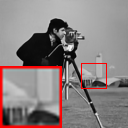}
		\small{34.91/0.951}
	\end{subfigure}%
	\hspace{1pt}
	\begin{subfigure}{0.15\linewidth}
		\centering
		\includegraphics[width=\linewidth]{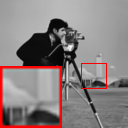}
		\small{35.04/0.954}
	\end{subfigure}
	
	\vspace{6pt}
	
	\begin{subfigure}{0.15\linewidth}
		\centering
		\includegraphics[width=\linewidth]{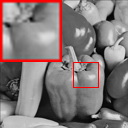}
		\small{PSNR/SSIM}
	\end{subfigure}%
	\hspace{1pt}
	\begin{subfigure}{0.15\linewidth}
		\centering
		\includegraphics[width=\linewidth]{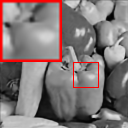}
		\small{33.27/0.949}
	\end{subfigure}%
	\hspace{1pt}
	\begin{subfigure}{0.15\linewidth}
		\centering
		\includegraphics[width=\linewidth]{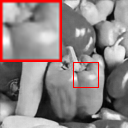}
		\small{34.23/0.955}
	\end{subfigure}%
	\hspace{1pt}
	\begin{subfigure}{0.15\linewidth}
		\centering
		\includegraphics[width=\linewidth]{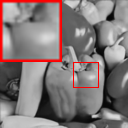}
		\small{34.44/0.961}
	\end{subfigure}%
	\hspace{1pt}
	\begin{subfigure}{0.15\linewidth}
		\centering
		\includegraphics[width=\linewidth]{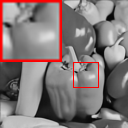}
		\small{34.40/0.960}
	\end{subfigure}%
	\hspace{1pt}
	\begin{subfigure}{0.15\linewidth}
		\centering
		\includegraphics[width=\linewidth]{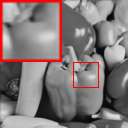}
		\small{34.55/0.961}
	\end{subfigure}
	
	\caption{Comparison of Image Restoration Effects of Different Methods at Noise Level $\alpha =27$.}
	\label{fig:phase_retrieval_comparison27}
\end{figure}
\begin{figure}[!t]
	\centering
	\captionsetup[subfigure]{skip=0.1pt} 
	\hspace{-\tabcolsep} 
	
	\begin{subfigure}{0.15\linewidth}
		\centering
		{Truth}
		\includegraphics[width=\linewidth]{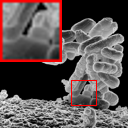}
		\small{PSNR/SSIM}
	\end{subfigure}%
	\hspace{1pt}
	\begin{subfigure}{0.15\linewidth}
		\centering
		{DeepMMSE}
		\includegraphics[width=\linewidth]{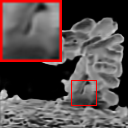}
		\small{23.49/0.850}
	\end{subfigure}%
	\hspace{1pt}
	\begin{subfigure}{0.15\linewidth}
		\centering
		{PrCom}
		\includegraphics[width=\linewidth]{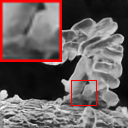}
		\small{25.21/0.786}
	\end{subfigure}%
	\hspace{1pt}
	\begin{subfigure}{0.15\linewidth}
		\centering
		{PRISTA-Net}
		\includegraphics[width=\linewidth]{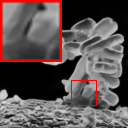}
		\small{25.55/0.866}
	\end{subfigure}%
	\hspace{1pt}
	\begin{subfigure}{0.15\linewidth}
		\centering
		{USP}
		\includegraphics[width=\linewidth]{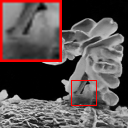}
		\small{25.28/0.884}
	\end{subfigure}%
	\hspace{1pt}
	\begin{subfigure}{0.15\linewidth}
		\centering
		{Ours}
		\includegraphics[width=\linewidth]{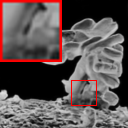}
		\small{25.68/0.902}
	\end{subfigure}
	
	\vspace{6pt}
	
	\begin{subfigure}{0.15\linewidth}
		\centering
		\includegraphics[width=\linewidth]{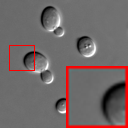}
		\small{PSNR/SSIM}
	\end{subfigure}%
	\hspace{1pt}
	\begin{subfigure}{0.15\linewidth}
		\centering
		\includegraphics[width=\linewidth]{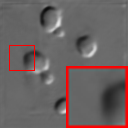}
		\small{27.52/0.900}
	\end{subfigure}%
	\hspace{1pt}
	\begin{subfigure}{0.15\linewidth}
		\centering
		\includegraphics[width=\linewidth]{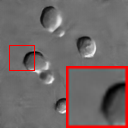}
		\small{29.41/0.922}
	\end{subfigure}%
	\hspace{1pt}
	\begin{subfigure}{0.15\linewidth}
		\centering
		\includegraphics[width=\linewidth]{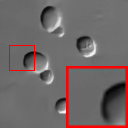}
		\small{33.12/0.954}
	\end{subfigure}%
	\hspace{1pt}
	\begin{subfigure}{0.15\linewidth}
		\centering
		\includegraphics[width=\linewidth]{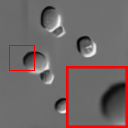}
		\small{33.63/0.960}
	\end{subfigure}%
	\hspace{1pt}
	\begin{subfigure}{0.15\linewidth}
		\centering
		\includegraphics[width=\linewidth]{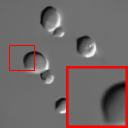}
		\small{33.75/0.958}
	\end{subfigure}
	
	\vspace{6pt}
	
	\begin{subfigure}{0.15\linewidth}
		\centering
		\includegraphics[width=\linewidth]{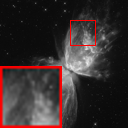}
		\small{PSNR/SSIM}
	\end{subfigure}%
	\hspace{1pt}
	\begin{subfigure}{0.15\linewidth}
		\centering
		\includegraphics[width=\linewidth]{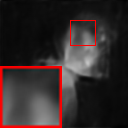}
		\small{27.70/0.792}
	\end{subfigure}%
	\hspace{1pt}
	\begin{subfigure}{0.15\linewidth}
		\centering
		\includegraphics[width=\linewidth]{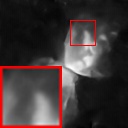}
		\small{29.10/0.843}
	\end{subfigure}%
	\hspace{1pt}
	\begin{subfigure}{0.15\linewidth}
		\centering
		\includegraphics[width=\linewidth]{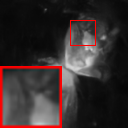}
		\small{30.38/0.837}
	\end{subfigure}%
	\hspace{1pt}
	\begin{subfigure}{0.15\linewidth}
		\centering
		\includegraphics[width=\linewidth]{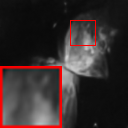}
		\small{30.20/0.841}
	\end{subfigure}%
	\hspace{1pt}
	\begin{subfigure}{0.15\linewidth}
		\centering
		\includegraphics[width=\linewidth]{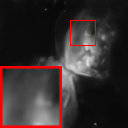}
		\small{30.68/0.853}
	\end{subfigure}
	
	\vspace{6pt}
	
	\begin{subfigure}{0.15\linewidth}
		\centering
		\includegraphics[width=\linewidth]{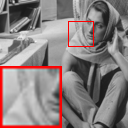}
		\small{PSNR/SSIM}
	\end{subfigure}%
	\hspace{1pt}
	\begin{subfigure}{0.15\linewidth}
		\centering
		\includegraphics[width=\linewidth]{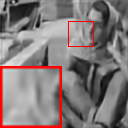}
		\small{25.52/0.833}
	\end{subfigure}%
	\hspace{1pt}
	\begin{subfigure}{0.15\linewidth}
		\centering
		\includegraphics[width=\linewidth]{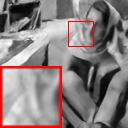}
		\small{26.67/0.868}
	\end{subfigure}%
	\hspace{1pt}
	\begin{subfigure}{0.15\linewidth}
		\centering
		\includegraphics[width=\linewidth]{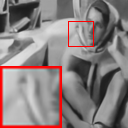}
		\small{27.69/0.875}
	\end{subfigure}%
	\hspace{1pt}
	\begin{subfigure}{0.15\linewidth}
		\centering
		\includegraphics[width=\linewidth]{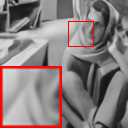}
		\small{27.73/0.881}
	\end{subfigure}%
	\hspace{1pt}
	\begin{subfigure}{0.15\linewidth}
		\centering
		\includegraphics[width=\linewidth]{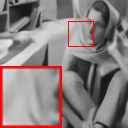}
		\small{28.05/0.878}
	\end{subfigure}
	
	\vspace{6pt}
	
	\begin{subfigure}{0.15\linewidth}
		\centering
		\includegraphics[width=\linewidth]{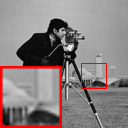}
		\small{PSNR/SSIM}
	\end{subfigure}%
	\hspace{1pt}
	\begin{subfigure}{0.15\linewidth}
		\centering
		\includegraphics[width=\linewidth]{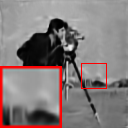}
		\small{26.10/0.815}
	\end{subfigure}%
	\hspace{1pt}
	\begin{subfigure}{0.15\linewidth}
		\centering
		\includegraphics[width=\linewidth]{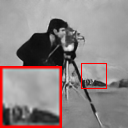}
		\small{27.75/0.875}
	\end{subfigure}%
	\hspace{1pt}
	\begin{subfigure}{0.15\linewidth}
		\centering
		\includegraphics[width=\linewidth]{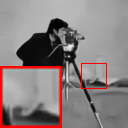}
		\small{28.36/0.875}
	\end{subfigure}%
	\hspace{1pt}
	\begin{subfigure}{0.15\linewidth}
		\centering
		\includegraphics[width=\linewidth]{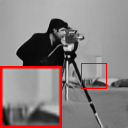}
		\small{28.75/0.874}
	\end{subfigure}%
	\hspace{1pt}
	\begin{subfigure}{0.15\linewidth}
		\centering
		\includegraphics[width=\linewidth]{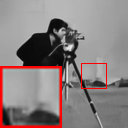}
		\small{28.95/0.878}
	\end{subfigure}
	
	\vspace{6pt}
	
	\begin{subfigure}{0.15\linewidth}
		\centering
		\includegraphics[width=\linewidth]{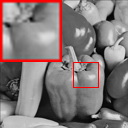}
		\small{PSNR/SSIM}
	\end{subfigure}%
	\hspace{1pt}
	\begin{subfigure}{0.15\linewidth}
		\centering
		\includegraphics[width=\linewidth]{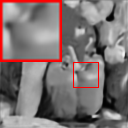}
		\small{25.61/0.834}
	\end{subfigure}%
	\hspace{1pt}
	\begin{subfigure}{0.15\linewidth}
		\centering
		\includegraphics[width=\linewidth]{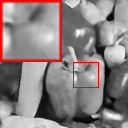}
		\small{26.79/0.871}
	\end{subfigure}%
	\hspace{1pt}
	\begin{subfigure}{0.15\linewidth}
		\centering
		\includegraphics[width=\linewidth]{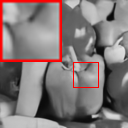}
		\small{27.80/0.870}
	\end{subfigure}%
	\hspace{1pt}
	\begin{subfigure}{0.15\linewidth}
		\centering
		\includegraphics[width=\linewidth]{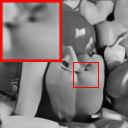}
		\small{27.99/0.878}
	\end{subfigure}%
	\hspace{1pt}
	\begin{subfigure}{0.15\linewidth}
		\centering
		\includegraphics[width=\linewidth]{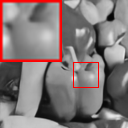}
		\small{28.10/0.880}
	\end{subfigure}
	
	\caption{Comparison of Image Restoration Effects of Different Methods at Noise Level $\alpha =81$ .}
	\label{fig:phase_retrieval_comparison}
\end{figure}
(1)~\textbf{Quantitative Evaluation}  
~The quantitative performance comparison of different methods under varying noise levels and with 4 masks on the 128×128 test dataset is detailed in Table  \ref{tab:pr_comparison}. We employed the average PSNR and SSIM across the 6 UNT-6 and 6 NT-6 images as evaluation metrics. The most effective result in each comparative row is highlighted in bold, with the second-highest performance underscored.

Experimental results confirm the proposed method’s superiority over four comparison algorithms across all noise levels and datasets. DLMMPR exhibits significant advantages on the 128×128 UNT-6 dataset, consistently outperforming others. Quantitatively, average PSNR on UNT-6 is   0.135 dB, 0.249 dB, and 0.134 dB
 higher than the second-best method at $\alpha = 9,27,81$, respectively. For NT-6, improvements are  0.153 dB, 0.179 dB, and 0.252dB. The increasing gain with noise intensity strongly validates DLMMPR’s robustness in high-noise scenarios.

PnP frameworks like DeepMMSE and PrComplex offer reconstruction improvements but suffer unstable performance in complex noise due to mismatches between prior models and actual noise distributions. USP, an unsupervised method, uses teacher-student distillation for noise-resistant training without ground truth, but its performance is highly sensitive to loss function-noise model adaptivity. Our comparison confirms our method achieves efficient computation and high reconstruction accuracy. Its advantages lie in superior numerical metrics, enhanced noise robustness, and better generalization.

Table \ref{tab:dataset_performance} data unequivocally shows the proposed method’s significant advantages under various noise levels and test images, stably leading in PSNR and SSIM. For the Yeast image in low-noise ($\alpha = 9$), our PSNR is  45.79
dB, 0.12  dB higher than USP (second-best). On ‘barbara’ and ‘peppers,’ PSNR gains further expand to between 0.13 and 0.30 dB. This indicates accurate detail preservation and high-quality reconstruction even with weak noise interference.

DLMMPR exhibits optimal performance with minimal fluctuations across different noise levels and images, demonstrating strong generalization ability. While comparison methods like USP achieve comparable results in some cases (e.g., SSIM on Yeast at 
 $\alpha = 9$), their overall stability is insufficient, especially on complex textures (‘Butterfly’) and high noise, where performance degrades significantly. These results fully verify DLMMPR’s effectiveness and superiority in image restoration, offering a more reliable solution for high-noise environments due to its dual advantages in accuracy and robustness.

(2)~\textbf{Qualitative Evaluation}  
~Figures \ref{fig:phase_retrieval_comparison9} - \ref{fig:phase_retrieval_comparison}
illustrate the visual comparison results. On the 128×128 test dataset with 4 masks and varying Poisson noise levels, our method preserves more structural details, while others exhibit over-smoothing or mottled backgrounds. These findings verify the significant value of the proposed matrix learning method for phase retrieval.

At $\alpha = 9$, DLMMPR significantly outperforms comparison methods on ‘boat’ and ‘cameraman’ images (PSNR: 40.04 dB and 41.60 dB; SSIM: 0.983 and 0.984 ). Visually, DLMMPR captures finer details (wood grain, mast ropes, beard, clothing wrinkles) with greater sharpness than blurred results from DeepMMSE and PrComplex. At $\alpha = 27$, DLMMPR reconstructs striped cloth periodicity and yeast cell details more accurately than PrComplex, avoiding adhesion blurring. At $\alpha = 81$, DLMMPR’s advantages are more pronounced: continuous bacterial cell contours and accurate intracellular light-dark distribution effectively suppress spot noise. DLMMPR’s performance gains over the second-best method increase from 0.25
dB (at $\alpha = 9$) to 0.4 dB (at $\alpha = 81$), confirming its robustness to high noise. This stems from accurately modeling image priors, enabling noise suppression via global constraints and detail preservation via local enhancement—achieving “denoising without damaging authenticity.”

In summation, the aforementioned experimental results unequivocally confirm that the proposed DLMMPR method demonstrates comprehensive leadership in both visual quality and quantitative metrics. This superior performance is consistently observed across a diverse range of image types and varying noise levels. These findings collectively serve to fully validate the effectiveness and inherent superiority of DLMMPR within the domain of phase retrieval tasks.
\section{Conclusions }

We propose a novel phase retrieval method that learns and optimizes the measurement matrix using target prior information via an end-to-end parameter learning framework. A fixed measurement matrix is optimized through learning, with all parameters converging to their optimal states. Extensive experiments verify the algorithm’s effectiveness. Future work will extend this learned optimization to other phase retrieval modules. For large-scale reconstruction, we will investigate lightweight networks and distributed optimization to reduce complexity while ensuring accuracy, enabling real-time imaging. These improvements aim to enhance generalization and practicality, supporting phase retrieval in precision measurement and remote sensing.

\section*{Data availability}
No data was used for the research described in the article.
\section*{Acknowledgments}
This research work was supported in part by the National Natural Science Foundation of China (grants 12261037).

\bibliographystyle{unsrt}
\bibliography{reference}

\end{document}